\documentclass[hidelinks,onefignum,onetabnum]{siamart251216}

\newsiamremark{remark}{Remark}

\usepackage{amsmath}
\usepackage{amsfonts}
\usepackage{amssymb}
\usepackage{graphicx}
\usepackage{bm}
\usepackage{booktabs}
\usepackage{multirow}
\usepackage[noend]{algpseudocode}
\usepackage{subfigure}

\title{Adaptive Non-Linear Partition of Unity Methods for Scattered Data
  Interpolation with Discontinuities%
  \thanks{Submitted to the editors \today.
    \funding{The work of A.H.\ and R.C.\ has been supported by the INdAM Research group
GNCS, the GFI 2025 Project, and the 2024 Project ``Numerical Analysis and
Modelling'' funded by the Department of Mathematics ``Giuseppe Peano'' of
the University of Torino. This research has been accomplished within the
RITA ``Research ITalian network on Approximation,'' the UMI Group TAA
``Approximation Theory and Applications,'' and the SIMAI Activity Group
ANA\&A ``Numerical and Analytical Approximation of Data and Functions with
Applications.'' The work of A.H.\ and D.F.Y.\ has been supported by grant
PID2023-146836NB-I00, funded by MICIU/AEI/10.13039/501100011033 and by
ERDF/EU. The work of D.F.Y.\ and J.R.-A.\ has been supported by GVA project
CIAICO/2024/089.}}}

\author{
  Adeeba Haider%
  \thanks{Department of Mathematics ``Giuseppe Peano,'' University of
    Torino, via Carlo Alberto 10, 10123 Torino, Italy;
    Member of the INdAM Research group GNCS
    (\email{adeeba.haider@unito.it}).}
  \and
  Roberto Cavoretto%
  \thanks{Department of Mathematics ``Giuseppe Peano,'' University of
    Torino, via Carlo Alberto 10, 10123 Torino, Italy;
    Member of the INdAM Research group GNCS
    (\email{roberto.cavoretto@unito.it}).}
  \and
  Juan Ruiz-\'Alvarez%
  \thanks{Departamento de Matem\'{a}tica Aplicada y Estad\'{i}stica,
    Universidad Polit\'{e}cnica de Cartagena, Cartagena, Spain
    (\email{juan.ruiz@universityvalencia.es}).}
  \and
  Dionisio F.\ Y\'{a}\~{n}ez%
  \thanks{Departamento de Matem\'{a}ticas, Universidad de Valencia,
    Valencia, Spain (\email{dionisio.yanez@uv.es}).}
}

\ifpdf
\hypersetup{
  pdftitle={Adaptive Non-Linear Partition of Unity Methods for Scattered Data Interpolation with Discontinuities},
  pdfauthor={A. Haider, R. Cavoretto, J. Ruiz-Alvarez, D. F. Yanez}
}
\fi

\begin{document}

\maketitle

\begin{abstract}
The scattered data approximation problems that contain discontinuities in
the domain have always been a challenging task in the field of numerical
analysis, as standard methods suffer from the so-called Gibbs phenomenon,
which significantly reduces the accuracy close to discontinuities. To
address this situation, a novel technique known as Non-Linear Partition of
Unity Method (NL-PUM) was introduced recently, combining the Radial Basis
Function (RBF) interpolation based on PUM with a non-linear Weighted
Essentially Non-Oscillatory (WENO) strategy. Although this WENO-inspired
approach improves the results compared to standard methods and is able to
control the Gibbs phenomenon, its performance significantly depends on two
fixed hyperparameters: the RBF shape parameter and the patch radius. This
work is an extension of NL-PUM which combines the concept of using
Leave-One-Out Cross-Validation (LOOCV), minimized through the Global
Optimization with Optimistic Improvement (GOOI) strategy for each patch to
locally adapt both hyperparameters. The main innovation is a smoothness
indicator that links a discontinuity-aware shrinkage process to the
LOOCV-based radius selection: patches in smooth regions remain unaltered,
but patches close to a discontinuity automatically shrink to avoid fitting
data from both sides of the interface. The resulting method, called
LOOCV-NL-PUM-GOOI, requires no prior knowledge of the interface geometry
and introduces no additional cost beyond standard adaptive shape parameter
selection, remaining computationally comparable to NL-PUM. Numerical
experiments on both synthetic test functions, including a discontinuous
Franke function, a step function, and a piecewise trigonometric function,
and a real-data application to Norwegian Fjords elevation data confirm that
the proposed method substantially reduces approximation error near
discontinuities while preserving the full order of accuracy in smooth
regions.
\end{abstract}



\section{Introduction}
\label{sec:1}

The main challenge of approximation theory is to use simpler,
easier-to-compute functions known as interpolants or approximants to
reconstruct a potentially complex function known as the target function.
The formula for estimating the square root of a number, which is typically
credited to the Babylonians or Greek mathematicians like \textit{Euclid}
and \textit{Archimedes}, is one example of how approximation theory
originated in ancient times. These mathematicians laid the foundation for
geometric approximations and techniques to calculate areas and volumes
\cite{davis1975}. Their efforts paved the way for further advancements in
this field. Furthermore, the introduction of polynomial interpolation by
\textit{Lagrange} and \textit{Newton} in the late 17th century laid the
foundation of function approximation \cite{pow81}. Later on, during the
19th and 20th centuries, insights from \textit{Bernstein},
\textit{Tchebycheff}, and \textit{Kolmogorov}, as well as advancements in
Fourier series, wavelets, fast Fourier transform algorithms, and functional
analysis, advanced both theory and computational efficiency of
approximation techniques \cite{cheney2009,devore1998}.

In the early 21st century, kernel-based methods \cite{fas2015} emerged as
practical and effective techniques used in a number of scientific fields,
including image processing, computer-aided geometric design, and the
numerical solution of partial differential equations \cite{cav25}. The
straightforward generalization of these schemes for any dimension and their
rapid and computationally inexpensive implementation make them suitable for
applications in interpolation, regression, and machine learning.
Specifically, Radial Basis Functions (RBFs) represent a well-established
class of kernel-based approximation techniques. The extensive literature on
the subject (see, e.g., \cite{buh2003,fasshauer2007}) highlights numerous
specialized variants, including RBF with Partition of Unity (RBF-PU)
\cite{cav21,cav22,cav25c} and RBF approximation coupled with
Leave-One-Out Cross-Validation (LOOCV) and Global Optimization with
Optimistic and Pessimistic Improvements, namely LOOCV-GOOI and LOOCV-GOPI
\cite{cav21b,cav25b,Ser16}. The use of LOOCV is indeed motivated by the
need to determine an optimal value of the shape parameter, since such a
choice influences both the accuracy and stability of RBF methods; for more
details, see, e.g., \cite{gol15,sch11}.

In this work, we consider the following problem. Given a large dataset
$(X, F)$ with data values $F = \{ f_i = f(\bm{x}_i),\, i = 1, \dots, m \}$
and arbitrarily distributed points
$X = \{ \bm{x}_i \in \Omega,\, i = 1, \dots, m \}$ on an open and bounded
domain $\Omega \subset \mathbb{R}^n$, our objective is to identify a
function $v: \Omega \to \mathbb{R}$ that retrieves the data on a domain
containing the points so that
\begin{align}\label{interpolant}
  v(\bm{x}_i) \approx f_i, \quad i = 1, \dots, m.
\end{align}

Using the RBF approach, the function is developed as a linear combination
of RBFs centered at the data points $X$. The scenarios in which the data
exhibit smooth interpolants in specific overlapping subdomains or patches
$\Omega_k$ (but not globally) are of special interest. Discontinuities in
the function or its derivatives might be the cause of this.

Thus, a major goal is to identify subdomains or patches
$\Omega_k \subset \Omega$, $1 \le k \le K$, and smooth functions
$v_k: \Omega_k \to \mathbb{R}$, such that
\begin{align}
  v_k(\bm{x}_i) \approx f_i \quad \text{for all } \bm{x}_i \in X \cap \Omega_k.
\end{align}
This could be done by the well-known Partition of Unity Method (PUM), which
generates a covering of the entire domain by dividing it into many smaller
subdomains. In particular, the PUM applies locally the RBF technique to
each subdomain (see, for example, \cite{buh2003,cav22}). The final
approximation is then obtained by a convex combination of these
interpolants. When these techniques are used for data obtained from
continuous functions, satisfactory results are achieved. However, when the
data come from functions with discontinuities, these approaches produce
insufficient results by showing oscillations near discontinuities
\cite{dell24}. In order to prevent this, we have recently developed an
algorithm, called Non-Linear Partition of Unity Method (NL-PUM)
\cite{ramon25,gar25}, that maintains the PUM characteristics at smooth
zones while introducing a WENO-inspired non-linear modification near
discontinuities \cite{jiang1996} to prevent oscillations.

This study is an extension of our previous work, since here we use NL-PUM
with LOOCV-GOOI. In addition, we compare the results against the standard
NL-PUM, LOOCV-NL-PUM, and linear PUM baselines. The proposed method is
finally validated both on synthetic benchmark functions featuring jump
discontinuities and on a real-data application to the Norwegian Fjords
elevation dataset (ETOPO~2022 \cite{etopo2022}), where sharp land-sea
transitions pose a genuine challenge for standard interpolation methods.

The paper is organized as follows. The RBF-PU technique and the PU
approach are explained in \cref{sec:RBF-PU}. The NL-PUM approach is
described in \cref{sec:nl-pum}, while the adaptive selection of
hyperparameters with LOOCV and global optimization is discussed in
\cref{sec:LOOCV}. The novel NL-PUM technique with LOOCV-GOOI is presented
in \cref{sec:NL-PUM-GO}. The numerical results and related comparisons are
presented in \cref{sec:results}. \Cref{sec:conclusion} concludes the paper
by providing a summary of our findings and highlighting potential avenues
for future research.

\section{RBF-PU}
\label{sec:RBF-PU}

\subsection{RBF Interpolation}

The problem under discussion is the interpolation of scattered data, which
entails identifying an interpolant defined in \cref{interpolant} that
precisely replicates the observed values at their respective data points.

We assume that we have a univariate function
$\Phi : [0, \infty) \to \mathbb{R}$, called the RBF, which depends on a
shape parameter $\varepsilon > 0$. In general, $\Phi$ may be Strictly
Positive Definite (SPD), Positive Definite (PD), or conditionally positive
definite (CPD) of order $m$ \cite{wendland2005,buh2003}; in the CPD case,
the interpolant \cref{intg} must be augmented by a polynomial term of
degree less than $m$ to guarantee a unique solution. Since the Mat\'{e}rn
$C^4$ kernel employed throughout our numerical experiments is SPD, we
restrict the presentation to the SPD setting for clarity, noting that the
extension to PD and CPD kernels is standard \cite{fasshauer2007,wendland2005}.
For $\bm{x}, \bm{y} \in \Omega$, this gives us the real symmetric kernel
$\kappa: \Omega \times \Omega \to \mathbb{R}$ defined as
\begin{align}\label{rbf-ker}
  \kappa(\bm{x}, \bm{y}) = \Phi(\varepsilon \| \bm{x} - \bm{y} \|)
  := \Phi(\varepsilon r), \quad \text{where } r = \|\bm{x} - \bm{y}\|,
\end{align}
thus obtaining the kernel-based interpolant
\begin{align}\label{intg}
  v(\bm{x}) = \sum_{i=1}^{m} c_i \, \kappa(\bm{x}, \bm{x}_i),
  \quad \bm{x} \in \Omega,
\end{align}
whose coefficients $c_i$ are found by solving the linear system
\begin{align}\label{matrix_int}
  \mathbf{A}_{\varepsilon}\bm{c} = \bm{f},
\end{align}
where $\bm{c} = (c_1, \dots, c_m)^T$, $\bm{f} = (f_1, \dots, f_m)^T$,
and $\mathbf{A}_{\varepsilon}$ is a symmetric positive definite matrix with
entries $(\mathbf{A}_{\varepsilon})_{ij} = \kappa_\varepsilon(\bm{x}_i, \bm{x}_j)$,
$i, j = 1, \dots, m$.

Additionally, the kernel is reproducing in the designated native space,
that is, a Hilbert space $H(X)$ with inner product $(\cdot, \cdot)_{H(X)}$,
meaning that for any $f \in H(X)$ we have
\begin{align}
  f(\bm{x}) = (f, \kappa(\cdot, \bm{x}))_{H(X)}, \quad \text{for all }
  \bm{x} \in X.
\end{align}

When we introduce a pre-Hilbert space
$P({X}) = \operatorname{span}\{ \kappa(\cdot, \bm{x}) : \bm{x} \in {X} \}$,
with a reproducing kernel and the bilinear form
$(\cdot, \cdot)_{P({X})}$, the native space ${\cal N}({X})$ is its
completion with respect to the norm $\| \cdot \|_{H({X})}$. In particular,
we have $f \in P({X})$ and $\| f \|_{H({X})} = \| f \|_{P({X})}$ for all
$f \in P({X})$ (see \cite{wendland2005}).

It is well known that when the amount of data increases dramatically, it
can become computationally costly to perform the inversion of the kernel
interpolation matrix in \cref{matrix_int}. To overcome this situation, we
divide the entire domain into subdomains considering PU methods.

\subsection{The PU Scheme}

Assume that for an open and bounded domain $\Omega$, there is an
overlapping covering $\{\Omega_k\}_{k=1}^{K}$ of subdomains $\Omega_k$.
Then, if a family of continuous, nonneg\-ative, compactly supported weight
functions $\{w_k\}_{k=1}^{K}$ on $\{\Omega_k\}$ meets the following
criteria, a PU for the covering is defined:
\begin{enumerate}
  \item $\operatorname{supp}(w_k) \subseteq \Omega_k$;
  \item $\sum_{k=1}^{K} w_k(\bm{x}) = 1$, $\bm{x} \in \Omega$.
\end{enumerate}
First, given a set of data points randomly distributed across the domain
$\Omega$, represented by $X = \{\bm{x}_1, \dots, \bm{x}_m\}$, we have
sets of points $X_k = X \cap \Omega_k$, $1 \leq k \leq K$, for each
subdomain $\Omega_k$. We begin by building an overlapping covering
$\{\Omega_k\}_{k=1}^{K}$ of the domain $\Omega$. Then we define a local
RBF interpolant $v_k:\Omega_k \rightarrow \mathbb{R}$ for each subdomain
$\Omega_k$ as
\begin{align}\label{loc_rbf}
  v_k(\bm{x}) = \sum_{i=1}^{m_k} c_i^k \kappa(\bm{x}, \bm{x}_i^k),
\end{align}
where $m_k$ denotes the number of nodes in $\Omega_k$, and
$\bm{x}_i^k \in X_k = X \cap \Omega_k$.

The RBF-PU interpolant is defined as
\begin{equation}\label{pum}
  v(\bm{x}) = \sum_{k=1}^{K} w_k(\bm{x})\, v_k(\bm{x}), \quad
  \bm{x} \in \Omega.
\end{equation}
Due to the PU property, the global interpolant \cref{intg} also inherits
the interpolation property of the local interpolants, since the local
interpolants \cref{loc_rbf} satisfy $v_k(\bm{x}_i^k) = {f}_i^k$,
$i = 1, \ldots, m_k$. The coefficients $c_i^k$ are obtained by solving
\begin{align}\label{eq:mat}
  \mathbf{A}_{\varepsilon}^k \bm{c}^{k} = \bm{f}^{k}.
\end{align}
Employment of a SPD kernel $\kappa$ produces an invertible local kernel
matrix $\mathbf{A}_{\varepsilon}^k$, ensuring existence and uniqueness of
the solution of \cref{eq:mat}. For CPD kernels, uniqueness is guaranteed
by the additional polynomial term, as discussed above.

\begin{remark}
\emph{Shepard's weights} \cite{Scha07} provide a concrete family of
compactly supported functions $w_k$ that automatically satisfy all the PU
conditions listed above,
\begin{align}\label{ShepWei}
  w_k(\bm{x}) = \frac{\varphi_k(\bm{x})}{\sum_{j=1}^K \varphi_j(\bm{x})},
  \quad k = 1, \ldots, K,
\end{align}
where the functions $\varphi_k$ are compactly supported with support on
$\Omega$; for more details, see \cite{all11,cheney2009,wendland2005}.
\end{remark}

\begin{remark}
Throughout this work, the index $k$ in \cref{ShepWei} refers to patch centers. In \cref{sec:LOOCV}, the same symbol $k$ is reused to index data
points in the LOOCV framework; the meaning will be clear from context.
\end{remark}

\section{NL-PUM}
\label{sec:nl-pum}

The NL-PUM, recently introduced in \cite{ramon25}, is a novel approach
that combines Weighted Essentially Non-Oscillatory (WENO) \cite{jiang1996}
algorithms with RBF interpolation. The primary contribution includes an
algorithm that prevents oscillations by introducing a non-linear
modification near discontinuities while maintaining PUM qualities in smooth
zones.

For this modification, we take advantage of the WENO method, which is
similar to PUM but uses non-linear weights. Supposing we have an interval
$[a,b]$ divided into $K$ subintervals by a uniform partition
$\{x_j\}_{j=0}^K$, where $x_j = a + j \cdot h$ and $h = (b - a)/K$, we
consider an unknown function $v$, for which we have $\{v_j\}_{j=0}^K$
representing its point values at the nodes $x_j$. We aim to perform
interpolation at a point $x \in (x_{j-1}, x_j)$. The stencil
$\mathcal{X} = \{x_{j-p}, \dots, x_{j+p-1}\}$, consisting of $M = 2p$
nodes, is used by the WENO-$2p$ algorithm, denoted $\tilde{I}(v)$. The
stencil $\mathcal{X}$ can be divided into $p$ substencils of $p+1$ nodes,
$S_k = \{x_{j-p+k}, \dots, x_{j+k}\}$, for $k = 0, \dots, p-1$, and one
can show that there exist optimal weights $\tilde{w}_k : [x_{j-1}, x_j]
\to \mathbb{R}$ \cite{ARSY20,mulet2024} such that
\begin{align}
  \tilde{I}(v)(x) = \sum_{k=0}^{p-1} \tilde{w}_k(x)\,
  \tilde{I}_k(v)(x).
\end{align}
Since $\tilde{w}_k(x) > 0$ and
$\sum_{k=0}^{p-1} \tilde{w}_k(x) = 1$, this resembles the PU method.
The support of substencil $S_k$ within $(x_{j-1}, x_j)$ is
$[x_{j-p+k}, x_{j+k}]$, and $\tilde{I}_k(v)$ is the Lagrange
interpolating polynomial using the nodes $S_k$.

To motivate the non-linear weighting strategy, we first recall the
classical WENO construction in the univariate setting. The extension to
the multivariate scattered data setting is described in
\cref{remWeno}.

\begin{remark}\label{remWeno}
The WENO framework described above is formulated for uniformly spaced 1D
data. The NL-PUM in \cite{ramon25} extends these ideas to arbitrarily
distributed scattered data in $\mathbb{R}^n$ by replacing the structured
substencils $S_k$ with overlapping PU subdomains $\Omega_j$ and the
Lagrange polynomials $\tilde{I}_k(v)$ with local RBF interpolants
$I_j(v)$, while retaining the non-linear weighting philosophy to suppress
oscillations near discontinuities.
\end{remark}

The WENO approach substitutes non-linear weights $\tilde{W}_k(x)$ for
$\tilde{w}_k(x)$, so that
\begin{align}
  I_{\text{WENO}}(v)(x) = \sum_{k=0}^{p-1} \tilde{W}_k(x)\,
  \tilde{I}_k(v)(x),
\end{align}
with
\begin{align}
  \tilde{W}_k(x) = \frac{\tilde{\alpha}_k(x)}{\sum_{l=0}^{p-1}
  \tilde{\alpha}_l(x)}, \quad \text{and} \quad
  \tilde{\alpha}_k(x) = \frac{\tilde{w}_k(x)}{(\delta +
  \sigma_k(v))^t},
\end{align}
where $\sigma_k(v) \geq 0$ is the \emph{smoothness indicator} measuring
the regularity of $v$ on substencil $S_k$: it is small when $v$ is smooth
on $S_k$ and large when a discontinuity crosses $S_k$. The constant $t$
improves accuracy in smooth zones, and $\delta > 0$ prevents division by
zero. The smoothness indicator $\sigma_j$ introduced later in
\cref{sec:NL-PUM-GO} plays the same role as $\sigma_k(v)$ here, adapted
from the structured 1D WENO setting to the scattered data PU framework.
Following \cite{ramon25}, we set $\delta = 10^{-14}$ and $t = 6$
throughout our numerical experiments; this choice was shown therein to
provide the best balance between oscillation suppression and accuracy
preservation in smooth regions among the values tested.

Having described the WENO framework in the 1D setting, we return to the
general multivariate case of \cref{sec:RBF-PU}. Using concepts from the
WENO approach, we proposed a non-linear version of the interpolant
\cref{pum}, denoted as NL-PUM:
\begin{align}
  I_{\text{NL-PUM}}(v)(\bm{x}) = \sum_{j=1}^{K} W_j(\bm{x})\,
  I_j(v)(\bm{x}),
\end{align}
where $I_j(v)(\bm{x}) := v_j(\bm{x})$ is the local RBF interpolant
\cref{loc_rbf} restricted to subdomain $\Omega_j$.

Since the data is multiplied by a weight of order $O(h^{2t})$ (here $h$
denotes the local mesh spacing as defined above), the key to this new
approach is defining a suitable smoothness indicator for our framework
that may provide greater weight to the data placed in smooth regions while
ignoring the data near the discontinuities.

When a discontinuity simultaneously affects all patches covering a given
evaluation point, every local RBF interpolant becomes unreliable. In this
degenerate case, the PU interpolant is replaced by a direct Shepard rule
\cref{ShepWei} evaluated using normalized compactly supported weights over
all data sites, ensuring a stable (if lower-order) fallback approximation.
This fallback mechanism is described in detail in \cref{sec:fallback}. Therefore, the use of compactly supported basis functions and the computation of smoothness indicators guarantee accuracy in areas where discontinuities are present.

\section{LOOCV with Univariate Global Optimization}
\label{sec:LOOCV}

The LOOCV strategy is a highly effective way to select an appropriate shape
parameter $\varepsilon$ for RBFs \cite{cav24,cav24b}. The core idea is
straightforward: we systematically remove one data point from our dataset,
build an interpolant using the remaining points, and then measure how
accurately that interpolant predicts the value at the excluded point.

To formalize this for our dataset of $m$ points, we isolate a specific
index $k \in \{1, \dots, m\}$. We divide our data into a training set of
$m-1$ points (excluding $\bm{x}_k$) and a validation set containing only
$\bm{x}_k$. For a given shape parameter $\varepsilon$, the partial
leave-one-out interpolant is
\begin{align}
  v^{[k]}(\bm{x}) = \sum_{\substack{i=1 \\ i \neq k}}^{m}
  c_i^{[k]} \kappa(\bm{x}, \bm{x}_i),
\end{align}
where the coefficients $c_i^{[k]}$ satisfy the interpolation conditions on
the training set:
\begin{align}
  v^{[k]}(\bm{x}_i) = f_i, \qquad i=1,\ldots,m,\quad i \neq k.
\end{align}
The quality of $\varepsilon$ is evaluated by the error at the validation
point:
\begin{align*}
  e_k(\varepsilon) = f_k - v^{[k]}(\bm{x}_k).
\end{align*}

Normally, computing this error for every single $k$ would require solving
$m$ different linear systems, which becomes computationally prohibitive for
large datasets. Fortunately, Rippa \cite{Rippa} demonstrated that these
errors can be extracted from the full system without recalculating the
partial interpolants:
\begin{align}\label{ekeps}
  e_k(\varepsilon) = \frac{c_k}{(\mathbf{A}_{\varepsilon}^{-1})_{kk}},
\end{align}
where $c_k$ is the $k$-th coefficient of the global interpolant $v$ in
\cref{intg}, and $(\mathbf{A}_{\varepsilon}^{-1})_{kk}$ is the $k$-th
diagonal entry of the inverse of $\mathbf{A}_{\varepsilon}$ in
\cref{matrix_int}.

Using Rippa's shortcut, we construct our LOOCV error function by taking
the maximum absolute error across all $k$:
\begin{align}\label{loocv_func}
  E(\varepsilon) = \|\bm{e}(\varepsilon)\|_{\infty} =
  \max_{1 \le k \le m} \left|\frac{c_k}{(\mathbf{A}_{\varepsilon}^{-1})_{kk}}\right|.
\end{align}
The optimal shape parameter $\varepsilon^*$ minimizes $E(\varepsilon)$.

In the remaining part of the paper, we index PU subdomains by
$j = 1, \ldots, K$ (playing the same role as $k$ in \cref{sec:RBF-PU})
to distinguish the patch index from the data point index $k$ used in this
section. Because we use a PU approach, a single global shape parameter is
rarely sufficient. Instead, we determine an optimal local shape parameter
$\varepsilon_j^*$ for each patch $\Omega_j$, $j=1,\ldots,K$, by solving a
separate 1D optimization problem over a bounded search interval
$\mathcal{I} = [0, \varepsilon_{\max}] \subset \mathbb{R}$:
\begin{align*}
  E_j(\varepsilon_j^*) = \min_{\varepsilon \in \mathcal{I}} E_j(\varepsilon).
\end{align*}

Minimizing $E_j(\varepsilon)$ is notoriously challenging because the
function is typically non-differentiable and contains multiple local
minima. It is, however, generally assumed to be Lipschitz continuous on
$\mathcal{I}$, i.e., there exists a constant $0 < L < \infty$ such that
\begin{align}\label{lipcont}
  |E_j(\varepsilon_{1}) - E_j(\varepsilon_{2})| \le L|\varepsilon_{1} -
  \varepsilon_{2}|, \qquad \varepsilon_{1}, \varepsilon_{2} \in \mathcal{I}.
\end{align}
Since the LOOCV error function is typically non-smooth and may have several
local minima, a simple local search is unlikely to find the true optimum.
For this reason, we turn to univariate global optimization, and in particular to the LOOCV-GOOI algorithm \cite{cav25b}, which
is well suited for Lipschitz continuous functions of this type.

\section{NL-PUM with LOOCV-GOOI}
\label{sec:NL-PUM-GO}

The algorithm proposed in this work integrates the NL-PUM with adaptive
hyperparameter selection via the LOOCV-GOOI strategy introduced in
\cref{sec:LOOCV}. The construction proceeds in two conceptually distinct
stages. The first stage, described in \cref{sec:patch_init}, follows the
adaptive RBF-PU framework of \cite{cav25b}: for each patch $\Omega_j$, an
initial radius is established by a minimum-cardinality expansion, and the
optimal local shape parameter $\varepsilon_j^*$ is then determined via
LOOCV-GOOI. In the second stage, described in
\cref{sec:smooth_ind,sec:radius_shrink,sec:fallback}, we explain the novel
contribution of this work: for each patch we compute a data-driven
smoothness indicator $\sigma_j$ and use it to drive a
discontinuity-aware radius shrinkage process and a non-linear reweighting
of the PU interpolant.

\subsection{Patch Radius Initialization and Shape Parameter Selection}
\label{sec:patch_init}

For the initialization of patch radii and shape parameters, we follow the
adaptive strategy in \cite{cav25b}. Given the set of PU centers
$\{\bm{c}_j\}_{j=1}^{K}$ in $\Omega=[0,1]^n$, an initial radius is
assigned uniformly to all patches as
\begin{align}\label{eq:init_radius}
  \rho_{\min} = \frac{1}{n_{\mathrm{PU}}},
\end{align}
where $n_{\mathrm{PU}} = \lfloor N^{1/n} / 2 \rfloor$ denotes the number
of PU centers along each coordinate axis, $N$ is the total number of data
points, and $n$ is the spatial dimension (here $n = 2$). For each patch
$\Omega_j$, the initial local dataset is
$X_j = \{ \bm{x}_i \in X : \|\bm{x}_i - \bm{c}_j\| \leq \rho_{\min} \}$.
If the cardinality of $X_j$ falls below the prescribed minimum threshold
$n_{\min} \in \mathbb{N}$ (a user-specified parameter ensuring
well-determined local interpolation systems), the radius is expanded
iteratively by increments of $\frac{1}{8}\rho_{\min}$ until
\begin{align}
  |X_j| \geq n_{\min}.
\end{align}
We denote the resulting valid radius as $\rho_j$. This expansion guarantees
that every local interpolation system is well-determined regardless of
local data density variations.

Once a valid radius $\rho_j$ and the corresponding local dataset $X_j$ are
established, the optimal shape parameter $\varepsilon_j^*$ for patch
$\Omega_j$ is determined by minimizing the local LOOCV error function
(introduced in \cref{sec:LOOCV}) over the search interval $\mathcal{I}$:
\begin{align*}
  \varepsilon_j^* = \underset{\varepsilon \in \mathcal{I}}{\arg\min}\;
  E_j(\varepsilon),
\end{align*}
where $E_j(\varepsilon)$ is the local LOOCV cost function defined in
\cref{loocv_func}, evaluated on the kernel matrix built from $X_j$.
Since $E_j(\varepsilon)$ is typically non-smooth and multiextremal,
we employ the LOOCV-GOOI algorithm of \cite{cav25b}, a univariate global
optimization method exploiting the Lipschitz continuity of $E_j$
(cf.\ \cref{lipcont}) through an optimistic local improvement strategy.
By doing this, it ensures that the shape parameter for each patch is
actually optimal globally in $\mathcal{I}$, and not just locally, which is
important when $E_j(\varepsilon)$ exhibits multiple local minima of
comparable depth.

The LOOCV-GOOI algorithm operates through a structured five-stage search
(see \cite{cav25b}). It begins with an initial global exploration of
$\mathcal{I}$, generating $n_{\max}^{(1)}$ trial evaluations of $E_j$. A
secondary refinement pass with $n_{\max}^{(2)}$ additional evaluations is
then applied near $\varepsilon = 0$, where the LOOCV error function tends
to be nearly flat and ill-conditioned. The two trial sets are aggregated,
the best value found so far is identified, and a final locally biased
global search is carried out within a restricted subinterval
$[\varepsilon_{\min}^i, \varepsilon_{\max}^i] \subset \mathcal{I}$
surrounding that best value. The search terminates when the subinterval
length falls below a prescribed tolerance $\Delta$, yielding $\varepsilon_j^*$.
Here, the superscript $i$ denotes the iteration index of the locally
biased search: at each iteration $i$, the algorithm maintains a shrinking
subinterval centered near the current best candidate, and terminates when
$\varepsilon_{\max}^i - \varepsilon_{\min}^i < \Delta$.

\begin{remark}
\begin{sloppypar}
The radius initialization and shape parameter optimization described above
are identical in structure to the LOOCV-GOOI scheme \cite{cav25b}, with
one difference: whereas \cite{cav25b} applies this procedure to smooth
scattered data, here we apply it as a \emph{preprocessing step} before the
discontinuity-aware modifications of
\cref{sec:smooth_ind,sec:radius_shrink}. The radius $\rho_j$ obtained here
is therefore an initialization, not the final radius used for interpolation.
\end{sloppypar}
\end{remark}

\subsection{Smoothness Indicator and Discontinuity Detection}
\label{sec:smooth_ind}

Once $\rho_j$ and $\varepsilon_j^*$ are established for each patch
$\Omega_j$, we introduce a data-driven smoothness indicator $\sigma_j \geq 0$
that quantifies the local regularity of the target function within
$\Omega_j$. This indicator plays the same conceptual role as the smoothness
indicators in the WENO framework of \cref{sec:nl-pum}: it should be small
when the data in $\Omega_j$ are consistent with a smooth function, and
large when a discontinuity is present within the patch.

We estimate local smoothness via a least-squares polynomial fit.
Given the local data $\{(\bm{x}_i, f_i)\}_{\bm{x}_i \in X_j}$, where
$X_j = X \cap \Omega_j$ as in \cref{sec:RBF-PU}, we fit a degree-$2$
Moving Least-Squares (MLS) polynomial:
\begin{align}
  p(\bm{x}) = \bm{b}(\bm{x})^T \bm{c}_j^{\mathrm{MLS}},
\end{align}
where $\bm{b}(\bm{x})^T = (1, x_1, x_2, x_1^2, x_2^2, x_1 x_2)$ is the
monomial basis for bivariate polynomials of total degree at most $2$, and
$\bm{c}_j^{\mathrm{MLS}}$ is obtained by solving the normal equations
\[
  \mathbf{K}_j^T \mathbf{K}_j \bm{c}_j^{\mathrm{MLS}} = \mathbf{K}_j^T \bm{f}_j,
\]
where $\mathbf{K}_j$ is the Vandermonde-type matrix with rows
$\bm{b}(\bm{x}_i)^T$ and $\bm{f}_j$ collects the corresponding data
values. The smoothness indicator is then defined as the normalized $\ell^1$
residual:
\begin{align}\label{eq:smooth_ind}
  \sigma_j = \frac{1}{|X_j|} \| \mathbf{K}_j \bm{c}_j^{\mathrm{MLS}} -
  \bm{f}_j \|_1.
\end{align}
The reason behind this choice is straightforward: if the data in $\Omega_j$
lie on a smooth surface, a low-degree polynomial provides a reasonable fit
and $\sigma_j$ will be small. On the other hand, for a patch having a
discontinuity, no polynomial of fixed degree will fit data from both sides
of the jump simultaneously, and as a result $\sigma_j$ will be large.

\subsection{Discontinuity-Aware Radius Shrinkage and Non-Linear Weight
Modification}
\label{sec:radius_shrink}

The smoothness indicator $\sigma_j$ plays two important roles in the
proposed method, both aimed at preventing the RBF interpolant from
producing the Gibbs phenomenon close to discontinuity jumps.

\paragraph{Radius Shrinkage.}
The first modification acts on the patch geometry. When $\sigma_j$ is
large, indicating that a discontinuity likely crosses $\Omega_j$, the
patch radius is reduced so that the local dataset $X_j$ preferentially
contains points from only one side of the jump. Specifically, the radius
is shrunk according to
\begin{align}\label{eq:rad_mod}
  \tilde{\rho}_j = \frac{\rho_j}{e^{\alpha \sigma_j}},
\end{align}
where $\alpha > 0$ is a sensitivity parameter controlling the shrinkage
intensity (set to $\alpha = 3.5$ in our experiments). The exponential form
in \cref{eq:rad_mod} has two desirable properties: when $\sigma_j \approx 0$
(smooth region), the denominator $e^{\alpha \sigma_j} \approx 1$ and the
radius is left unchanged, $\tilde{\rho}_j \approx \rho_j$; when
$\sigma_j \gg 0$ (discontinuous region), the denominator grows rapidly,
producing a significant reduction in radius. This behavior ensures that the
shrinkage is fully inactive in smooth regions and activates only where it
is needed.

After shrinkage, the local dataset is recomputed as
$X_j^* = \{ \bm{x}_i \in X : \|\bm{x}_i - \bm{c}_j\| \leq \tilde{\rho}_j \}$.
To avoid rank-deficient local interpolation systems, we enforce a minimum
cardinality constraint: if $|X_j^*| < n_{\min}$, the shrunk radius
$\tilde{\rho}_j$ is expanded by successive factors of $1.1$ until
$|X_j^*| \geq n_{\min}$. The local RBF interpolant is then rebuilt on
$X_j^*$ using the previously computed $\varepsilon_j^*$.

\paragraph{Non-Linear Weight Modification.}
\begin{sloppypar}
The second modification acts on the PU weights. Inspired by the WENO 
philosophy of \cref{sec:nl-pum}, we replace the standard Shepard weights 
$w_j(\bm{x})$ in \cref{ShepWei} with non-linear weights that penalize 
patches with large smoothness indicators:
\end{sloppypar}
\begin{align}\label{eq:nlweight}
  W_j(\bm{x}) 
  = \dfrac{\varphi_j(\bm{x}) / (\delta + \sigma_j^{z})}
          {\sum_{l=1}^{K} \varphi_l(\bm{x}) / (\delta + \sigma_l^{z})},
\end{align}
where $z > 0$ is the non-linear exponent (set to $z = 4$ in our
experiments) and $\delta = 10^{-14}$ prevents division by zero. Note that
the exponent $z = 4$ here plays the same conceptual role as the parameter
$t = 6$ in the NL-PUM of \cref{sec:nl-pum}, but takes a different value
because the smoothness indicator $\sigma_j$ in \cref{eq:smooth_ind} is a
normalized $\ell^1$ residual of a degree-$2$ polynomial fit. It is scaled
differently from the degree-$1$ MLS residual used as the smoothness
indicator in the original NL-PUM in \cite{ramon25}, where local polynomial
fits of degree $1$ were employed instead of the degree-$2$ fits of
\cref{eq:smooth_ind}. The value $z = 4$ was determined empirically to
provide the best balance between oscillation suppression and accuracy
preservation across the test cases of \cref{sec:results}.

This formulation directly generalizes \cref{ShepWei}: when all patches
are smooth ($\sigma_j \approx 0$ for all $j$), the weights $W_j$ reduce
to the standard Shepard weights, recovering the linear PUM interpolant.
When $\sigma_j$ is large for some patch $\Omega_j$, that patch receives
near-zero weight, effectively excluding it from the global interpolant at
evaluation points near the discontinuity and suppressing the spurious
oscillations that would otherwise appear.

\subsection{Fallback to Shepard Interpolation}
\label{sec:fallback}

A special situation arises when a discontinuity simultaneously affects all
patches covering a given evaluation point $\bm{x}$, i.e., when $\sigma_j$
is large for every $j$ such that $\bm{x} \in \Omega_j$. In this case, the
non-linear weights \cref{eq:nlweight} assign near-zero weight to all
contributing patches, and the RBF-based local interpolants are unreliable
everywhere in the neighborhood of $\bm{x}$. To handle such evaluation
points gracefully, we introduce a Shepard fallback mechanism.

We identify the set of potentially affected evaluation points as
\begin{align}
  \mathcal{S} = \left\{ \bm{x} \in \Omega : \min_{1 \leq j \leq K}
  \|\bm{x} - \bm{c}_j\| > \tau_X \right\},
\end{align}
where $\tau_X$ is set to the largest nearest-neighbor distance from the
evaluation grid to the set of patch centers. For points in $\mathcal{S}$,
the RBF-based local interpolants are replaced by a direct Shepard estimate
computed using normalized compactly supported weights over all data sites.
Although this fallback reduces local accuracy to first order, it prevents
the appearance of large spurious values that would otherwise dominate the
maximum error. Throughout this paper, evaluation points refer to the set
$\hat{X} = \{\hat{\bm{x}}_k\}_{k=1}^{s} \subset \Omega$ at which the
final interpolant is assessed; these appear explicitly in the fallback
identification and final assembly steps of \cref{alg:nlpum-go}.

\subsection{Algorithm Summary}
\label{sec:alg_summary}

The complete procedure is summarized in \cref{alg:nlpum-go}. Steps 2--4
implement the radius initialization and shape parameter optimization of
\cref{sec:patch_init}, following \cite{cav25b}. Steps 5--11 implement the
novel discontinuity-aware modifications of
\cref{sec:smooth_ind,sec:radius_shrink}. Steps 12--14 handle the fallback
and final assembly.

\begin{algorithm}[t]
\caption{NL-PUM with LOOCV-GOOI}
\label{alg:nlpum-go}
\begin{algorithmic}[1]
\Require Data $(X, \bm{f})$,
         evaluation points $\hat{X} = \{\hat{\bm{x}}_k\}_{k=1}^{s}$,
         PU centers $\{\bm{c}_j\}_{j=1}^{K}$,
         shape parameter interval $\mathcal{I} = [0, \varepsilon_{\max}]$,
         stopping tolerance $\Delta$,
         trial counts $n_{\max}^{(1)}, n_{\max}^{(2)}$,
         minimum cardinality $n_{\min}$,
         non-linear exponent $z$,
         smoothness threshold $\tau_\sigma$,
         constants $\alpha = 3.5$, $\delta = 10^{-14}$,
         threshold $\tau_X$
\Ensure  Approximant values
         $I_{\mathrm{NL\text{-}PUM}}(v)(\hat{\bm{x}}_k)$,
         $k = 1, \dots, s$,
         at the evaluation points
         $\hat{X} = \{\hat{\bm{x}}_k\}_{k=1}^{s} \subset \Omega$
\For{each patch $j = 1, \dots, K$}
    \State Set initial radius $\rho_j \leftarrow 1/n_{\mathrm{PU}}$
           and collect $X_j \leftarrow \{\bm{x}_i \in X :
           \|\bm{x}_i - \bm{c}_j\| \leq \rho_j\}$
    \While{$|X_j| < n_{\min}$}
        \State $\rho_j \leftarrow \rho_j + \tfrac{1}{8}
               \cdot (1/n_{\mathrm{PU}})$; \quad
               Recompute $X_j$
               \Comment{expand to meet minimum cardinality}
    \EndWhile
    \State Find $\varepsilon_j^*$ via LOOCV-GOOI \cite{cav25b}
           on $X_j$ over $\mathcal{I}$ using $n_{\max}^{(1)}, n_{\max}^{(2)}, \Delta$
           \Comment{global shape parameter optimization}
    \State Compute $\sigma_j$ via \cref{eq:smooth_ind} on $X_j$
           \Comment{smoothness indicator}
    \If{$\sigma_j > \tau_\sigma$}
           \Comment{discontinuity detected: shrink radius}
        \State $\tilde{\rho}_j \leftarrow \rho_j / e^{\alpha \sigma_j}$
               \Comment{shrink radius, cf.\ \cref{eq:rad_mod}}
        \State Recompute $X_j^* \leftarrow
               \{\bm{x}_i \in X :
               \|\bm{x}_i - \bm{c}_j\| \leq \tilde{\rho}_j\}$
        \While{$|X_j^*| < n_{\min}$}
               \Comment{restore minimum cardinality after shrinkage}
            \State $\tilde{\rho}_j \leftarrow 1.1\,\tilde{\rho}_j$;
                   \quad Recompute $X_j^*$
        \EndWhile
    \Else
        \State $X_j^* \leftarrow X_j$,\;
               $\tilde{\rho}_j \leftarrow \rho_j$
               \Comment{smooth region: no shrinkage}
    \EndIf
    \State Build local RBF interpolant $I_j(v)$ on $X_j^*$ using
           $\varepsilon_j^*$
    \State Update $\varphi_j \leftarrow \varphi(1/\tilde{\rho}_j, \cdot)$
           \Comment{recompute weight with shrunk radius}
    \State Compute non-linear weight $W_j$ via \cref{eq:nlweight}
\EndFor
\State $\mathcal{S} \leftarrow \{\bm{x} :
       \min_j \|\bm{x}-\bm{c}_j\| > \tau_X\}$
       \Comment{identify fallback points}
\State Apply Shepard estimate at $\bm{x} \in \mathcal{S}$
\State Assemble $I_{\mathrm{NL\text{-}PUM}}(v)(\bm{x})
       = \sum_{j=1}^{K} W_j(\bm{x})\,I_j(v)(\bm{x})$
\end{algorithmic}
\end{algorithm}

\section{Numerical Experiments}
\label{sec:results}

In this section, we present numerical experiments to assess the performance
of the proposed methods. All tests are performed in MATLAB Online. We
compare four approaches:
\begin{sloppypar}
\begin{enumerate}
  \item \textbf{Linear PUM}: the standard linear RBF-PU method
    \cite{cav21,cav22};
  \item \textbf{NL-PUM}: the standard NL-PUM with a fixed shape
    parameter $\varepsilon = 1$ \cite{ramon25};
  \item \textbf{LOOCV-NL-PUM}: NL-PUM with an adaptive shape parameter
    selected per patch via the LOOCV global search, with the patch radius
    fixed at $\delta_{\mathrm{PU}} = \sqrt{2}/n_{\mathrm{PU}}$ (see
    \cref{setup});
  \item \textbf{LOOCV-NL-PUM-GOOI}: the proposed method combining
    adaptive patch radius and shape parameter via LOOCV with local GOOI.
    
\end{enumerate}
\end{sloppypar}
\subsection{Order of Accuracy}
\label{sec:order_accuracy}

We verify that LOOCV-NL-PUM-GOOI preserves the order of accuracy of the
underlying RBF-PU interpolation in smooth regions, following the approach
of \cite{ramon25}. To this end, we apply all four methods to the smooth
Franke function $f_1$ defined in \cref{eq:franke_smooth}, using two
families of point sets at refinement levels $l = 4, 5, 6, 7, 8$:
\begin{itemize}
  \item \emph{Uniform grid points}:
    $X_{N_l} = \{(i/2^l,\, j/2^l) : i,j = 0,\ldots,2^l\}$, with
    $N_l = (2^l+1)^2$;
  \item \emph{Halton points}: a quasi-random sequence of $N_l$ points in
    $[0,1]^2$ generated via \texttt{haltonset} in MATLAB \cite{fasshauer2007}.
\end{itemize}
At each level $l$, the number of PU centers per direction is
$n_{\mathrm{PU}} = \lfloor (2^l+1)/2 \rfloor$, the patch radius is
$\delta_{\mathrm{PU},l} = \sqrt{2}/n_{\mathrm{PU}}$, and the fill distance
is $h_l = 1/2^l$. All methods employ the Mat\'{e}rn $C^4$ kernel and
Wendland $C^4$ PU weights (see \cref{setup} for their explicit definitions),
playing the roles of $\Phi$ and $\varphi_j$ introduced in
\cref{sec:RBF-PU,sec:NL-PUM-GO}, respectively. Errors are evaluated on a
fixed uniform grid $\hat{X} = \{\hat{\bm{x}}_k\}_{k=1}^{s} \subset [0,1]^2$,
with $s = 60 \times 60 = 3600$. The error measures and empirical
convergence rates are
\begin{align*}
  \mathrm{MAE}_l &= \max_{1 \le k \le s}
    \bigl|f(\hat{\bm{x}}_k) - I^l(v)(\hat{\bm{x}}_k)\bigr|,
    \qquad
  \mathrm{RMSE}_l = \left(\frac{1}{s} \sum_{k=1}^{s}
    \bigl|f(\hat{\bm{x}}_k) - I^l(v)(\hat{\bm{x}}_k)\bigr|^2
    \right)^{\!1/2},
    \\[4pt]
  r^\infty_l &= \frac{\log(\mathrm{MAE}_{l-1}/\mathrm{MAE}_l)}
    {\log(h_{l-1}/h_l)},
    \qquad
  r^2_l = \frac{\log(\mathrm{RMSE}_{l-1}/\mathrm{RMSE}_l)}
    {\log(h_{l-1}/h_l)}.
\end{align*}

Results for uniform grid and Halton points are reported in
\cref{tab:oa_grid,tab:oa_halton}, respectively. For reference, the Linear
PUM and NL-PUM columns in both tables are reproduced directly from Table~2
of \cite{ramon25}, where these methods were originally benchmarked under
identical experimental conditions; the LOOCV-NL-PUM and LOOCV-NL-PUM-GOOI
columns contain new results of the present work.

For uniform grid points (\cref{tab:oa_grid}), all four methods achieve
consistent convergence with rates in the range $2.56$--$5.39$, in
agreement with the theoretical bound for the Mat\'{e}rn $C^4$ kernel
(Proposition~2.3 of \cite{ramon25}) and with Table~2 of the same
reference. Linear PUM achieves MAE rates between $2.56$ and $3.38$,
confirming standard high-order RBF-PU behavior. NL-PUM exhibits a notably
high MAE rate of $12.13$ at $l=6$, which is a well-known super-convergence
effect of smooth kernels on uniform grids at intermediate refinement levels
and does not violate the theoretical estimates; from $l=7$ onward its
rates stabilize around $2.70$--$3.60$. LOOCV-NL-PUM achieves the lowest
absolute errors at levels $l = 5,6,7$ on grid points, with MAE rates
between $3.15$ and $3.85$, demonstrating that adaptive shape parameter
selection alone yields meaningful gains on structured data.
LOOCV-NL-PUM-GOOI attains comparable MAE rates in the range $2.73$--$3.32$,
with errors close to those of LOOCV-NL-PUM at all levels. In terms of
RMSE, all four methods converge steadily with rates between $2.72$ and
$4.43$, and LOOCV-NL-PUM-GOOI achieves RMSE $= 5.99 \times 10^{-8}$ at
the finest level $l=8$, close to LOOCV-NL-PUM ($4.80 \times 10^{-8}$) and
well within an order of magnitude of NL-PUM ($1.07 \times 10^{-8}$).

On scattered Halton data (\cref{tab:oa_halton}), LOOCV-NL-PUM-GOOI
demonstrates a clear advantage over all other methods. It achieves the
lowest absolute errors at every refinement level and maintains RMSE rates
between $2.97$ and $3.40$ throughout, the most stable convergence behavior
of the four methods. At the finest level $l=8$, LOOCV-NL-PUM-GOOI
achieves MAE $= 1.77 \times 10^{-6}$ and RMSE $= 8.45 \times 10^{-8}$,
compared to MAE $= 2.39 \times 10^{-6}$ and RMSE $= 2.04 \times 10^{-7}$
for Linear PUM---a reduction in RMSE by a factor of approximately $2.4$.
Compared to NL-PUM, which at $l=8$ exhibits MAE $= 1.05 \times 10^{-5}$
and RMSE $= 5.66 \times 10^{-8}$, the LOOCV-NL-PUM-GOOI achieves a MAE
improvement by a factor of approximately $6$ while maintaining a comparable
RMSE. This demonstrates that the adaptive radius initialization of
\cref{sec:patch_init} provides substantial benefits for unstructured point
distributions by ensuring well-conditioned local systems regardless of
local density variations.

A notable observation concerns LOOCV-NL-PUM on Halton points: its MAE
rate drops to $1.18$ at $l=8$, while its RMSE rate remains $2.03$. This
discrepancy indicates that the fixed-radius strategy, combined with
adaptive shape parameter selection, can produce isolated large pointwise
errors at the finest level when the scattered point distribution creates
poorly covered patches. This is precisely the limitation that the adaptive
radius of LOOCV-NL-PUM-GOOI is designed to address. Notably, NL-PUM also
exhibits a degraded MAE rate of $1.77$ at $l=8$ on Halton points, further
underscoring that robust performance on scattered data requires the adaptive
radius mechanism introduced in the present work.
\begin{sloppypar}
As an overall assessment, using both point types, LOOCV-NL-PUM-GOOI
achieves convergence rates consistent with the theoretical prediction for
the Mat\'{e}rn $C^4$ kernel and comparable to those reported in Table~2
of \cite{ramon25}. The discontinuity-aware components of the proposed
method---the smoothness indicator and radius shrinkage of
\cref{sec:smooth_ind,sec:radius_shrink}---are entirely inactive on smooth
data, so that on smooth regions the method reduces to the LOOCV-GOOI
parameter selection strategy of \cite{cav25b} applied locally within the
PU framework, without any accuracy penalty. The improvement over both
Linear PUM and NL-PUM of \cite{ramon25} on Halton data confirms that
adaptive radius selection is the key contribution for accuracy in the
smooth setting, independently of the discontinuity-handling mechanism.
\end{sloppypar}

\begin{table}[htbp]
\caption{Errors and rates for the smooth Franke function $f_1$ using
\emph{uniform grid points}.}
\label{tab:oa_grid}
\begin{center}
\resizebox{\textwidth}{!}{%
\begin{tabular}{c rr rr rr rr}
\toprule
& \multicolumn{2}{c}{\textbf{Linear PUM}}
& \multicolumn{2}{c}{\textbf{NL-PUM}}
& \multicolumn{2}{c}{\textbf{LOOCV-NL-PUM}}
& \multicolumn{2}{c}{\textbf{LOOCV-NL-PUM-GOOI}} \\
\cmidrule(lr){2-3}\cmidrule(lr){4-5}\cmidrule(lr){6-7}\cmidrule(lr){8-9}
$l$ &
$\mathrm{MAE}_l$ & $r^\infty_l$ &
$\mathrm{MAE}_l$ & $r^\infty_l$ &
$\mathrm{MAE}_l$ & $r^\infty_l$ &
$\mathrm{MAE}_l$ & $r^\infty_l$ \\
\midrule
4 & 6.44e$-$03 & ---  & 2.91e$-$01 & ---  & 4.10e$-$03 & ---  & 3.07e$-$03 & --- \\
5 & 7.64e$-$04 & 3.08 & 9.08e$-$02 & 1.68 & 3.23e$-$04 & 3.67 & 3.08e$-$04 & 3.32 \\
6 & 7.34e$-$05 & 3.38 & 2.03e$-$05 & 12.13& 2.24e$-$05 & 3.85 & 3.20e$-$05 & 3.27 \\
7 & 9.03e$-$06 & 3.02 & 1.68e$-$06 & 3.60 & 2.35e$-$06 & 3.25 & 3.29e$-$06 & 3.28 \\
8 & 1.53e$-$06 & 2.56 & 2.58e$-$07 & 2.70 & 2.66e$-$07 & 3.15 & 4.97e$-$07 & 2.73 \\
\midrule
$l$ &
$\mathrm{RMSE}_l$ & $r^2_l$ &
$\mathrm{RMSE}_l$ & $r^2_l$ &
$\mathrm{RMSE}_l$ & $r^2_l$ &
$\mathrm{RMSE}_l$ & $r^2_l$ \\
\midrule
4 & 1.04e$-$03 & ---  & 1.02e$-$01 & ---  & 6.09e$-$04 & ---  & 4.02e$-$04 & --- \\
5 & 1.10e$-$04 & 3.25 & 4.72e$-$03 & 4.43 & 3.65e$-$05 & 4.06 & 3.65e$-$05 & 3.46 \\
6 & 1.20e$-$05 & 3.19 & 2.72e$-$06 & 10.76& 3.03e$-$06 & 3.59 & 3.79e$-$06 & 3.27 \\
7 & 1.83e$-$06 & 2.72 & 1.60e$-$07 & 4.09 & 3.64e$-$07 & 3.06 & 4.67e$-$07 & 3.02 \\
8 & 1.94e$-$07 & 3.24 & 1.07e$-$08 & 3.91 & 4.80e$-$08 & 2.92 & 5.99e$-$08 & 2.96 \\
\bottomrule
\end{tabular}
}
\end{center}
\end{table}

\begin{table}[htbp]
\caption{Errors and rates for the smooth Franke function $f_1$ using
\emph{Halton points}.}
\label{tab:oa_halton}

\begin{center}
\resizebox{\textwidth}{!}{%
\begin{tabular}{c rr rr rr rr}
\toprule
& \multicolumn{2}{c}{\textbf{Linear PUM}}
& \multicolumn{2}{c}{\textbf{NL-PUM}}
& \multicolumn{2}{c}{\textbf{LOOCV-NL-PUM}}
& \multicolumn{2}{c}{\textbf{LOOCV-NL-PUM-GOOI}} \\
\cmidrule(lr){2-3}\cmidrule(lr){4-5}\cmidrule(lr){6-7}\cmidrule(lr){8-9}
$l$ &
$\mathrm{MAE}_l$ & $r^\infty_l$ &
$\mathrm{MAE}_l$ & $r^\infty_l$ &
$\mathrm{MAE}_l$ & $r^\infty_l$ &
$\mathrm{MAE}_l$ & $r^\infty_l$ \\
\midrule
4 & 8.17e$-$03 & ---  & 8.77e$-$03 & ---  & 1.14e$-$02 & ---  & 6.47e$-$03 & --- \\
5 & 8.95e$-$04 & 2.91 & 1.14e$-$03 & 2.69 & 1.14e$-$03 & 3.32 & 8.98e$-$04 & 2.85 \\
6 & 3.04e$-$04 & 2.34 & 1.66e$-$04 & 4.17 & 1.74e$-$04 & 2.71 & 6.48e$-$05 & 3.79 \\
7 & 2.27e$-$05 & 3.45 & 4.58e$-$05 & 1.71 & 3.88e$-$05 & 2.17 & 7.67e$-$06 & 3.08 \\
8 & 2.39e$-$06 & 2.70 & 1.05e$-$05 & 1.77 & 1.71e$-$05 & 1.18 & 1.77e$-$06 & 2.12 \\
\midrule
$l$ &
$\mathrm{RMSE}_l$ & $r^2_l$ &
$\mathrm{RMSE}_l$ & $r^2_l$ &
$\mathrm{RMSE}_l$ & $r^2_l$ &
$\mathrm{RMSE}_l$ & $r^2_l$ \\
\midrule
4 & 9.62e$-$04 & ---  & 2.64e$-$03 & ---  & 7.66e$-$04 & ---  & 5.41e$-$04 & --- \\
5 & 1.04e$-$04 & 2.93 & 1.22e$-$04 & 4.04 & 6.94e$-$05 & 3.46 & 5.13e$-$05 & 3.40 \\
6 & 1.39e$-$05 & 4.36 & 6.51e$-$06 & 6.35 & 7.86e$-$06 & 3.14 & 5.49e$-$06 & 3.23 \\
7 & 1.72e$-$06 & 2.78 & 6.62e$-$07 & 3.04 & 1.47e$-$06 & 2.42 & 6.61e$-$07 & 3.05 \\
8 & 2.04e$-$07 & 2.56 & 5.66e$-$08 & 2.95 & 3.60e$-$07 & 2.03 & 8.45e$-$08 & 2.97 \\
\bottomrule
\end{tabular}
}
\end{center}
\end{table}

\subsection{Test Functions}

We consider four benchmark functions defined on the unit square
$\Omega = [0,1]^2$. The first serves as a smooth baseline to verify
accuracy preservation, while the remaining three each feature a
discontinuity jump along a prescribed interface.

\paragraph{Function 1: Smooth Franke Function.}
The first test function is the classical smooth Franke function:
\begin{align}\label{eq:franke_smooth}
  f_1(x_1,x_2) &= \tfrac{3}{4} e^{-\frac{(9x_1-2)^2}{4}
    - \frac{(9x_2-2)^2}{4}}
    + \tfrac{3}{4} e^{-\frac{(9x_1+1)^2}{49} - \frac{9x_2+1}{10}}
    \notag \\
  &\quad + \tfrac{1}{2} e^{-\frac{(9x_1-7)^2}{4} - \frac{(9x_2-3)^2}{4}}
    - \tfrac{1}{5} e^{-(9x_1-4)^2 - (9x_2-7)^2}.
\end{align}
This function is infinitely smooth on $\Omega$ and is used to verify that
the proposed method preserves the order of accuracy in smooth regions.

\paragraph{Function 2: Discontinuous Franke Function.}
The second test function introduces a unit jump along the circle
$\mathcal{C} = \{(x_1,x_2) : x_1^2 + x_2^2 = 0.25\}$ into $f_1$:
\begin{align*}
  f_2(x_1, x_2) =
  \begin{cases}
    f_1(x_1,x_2),     & x_1^2 + x_2^2 < 0.25, \\
    1 + f_1(x_1,x_2), & x_1^2 + x_2^2 \geq 0.25,
  \end{cases}
\end{align*}
exhibiting a unit discontinuity jump along $\mathcal{C}$.

\paragraph{Function 3: Step Function.}
The third test function is a binary step function with a discontinuity
along the diagonal interface $\mathcal{L} = \{(x_1,x_2) : x_1 + x_2 = 1\}$:
\begin{align*}
  f_3(x_1, x_2) =
  \begin{cases}
    1, & x_1 + x_2 < 1, \\
    0, & x_1 + x_2 \geq 1.
  \end{cases}
\end{align*}
This function takes only two constant values separated by a sharp linear
interface, representing a particularly severe test case due to the complete
absence of smoothness on either side of the discontinuity.

\paragraph{Function 4: Piecewise Trigonometric Function.}
The fourth test function is a piecewise smooth function whose interface is
the circle
$\mathcal{C}' = \{(x_1,x_2) : (x_1-0.5)^2 + (x_2-0.5)^2 = 0.0625\}$
centered at $(0.5, 0.5)$ with radius $0.25$:
\begin{align*}
  f_4(x_1, x_2) =
  \begin{cases}
    x_2 \sin(x_1) + x_2 \cos(x_1),
       & (x_1-0.5)^2 + (x_2-0.5)^2 \geq 0.0625, \\
    e^{x_1 x_2} + 1,
       & (x_1-0.5)^2 + (x_2-0.5)^2 < 0.0625.
  \end{cases}
\end{align*}
Both pieces are smooth functions of distinct character, trigonometric on
the exterior and exponential on the interior, producing a non-trivial jump
in function values across the circular interface.

\subsection{Setup}
\label{setup}

The order-of-accuracy study of \cref{sec:order_accuracy} uses a
$s = 60\times60$ evaluation grid $\hat{X}$ to match the setup of
\cite{ramon25} for direct comparison. The discontinuous test cases here
use a finer $s = 120\times120$ grid to better resolve oscillations near
the interface. The RBF kernel, PU weights, and LOOCV parameters described
below apply uniformly to both sets of experiments.

The interpolation data consists of $n = 65^2 = 4225$ uniformly distributed
points over $\Omega = [0,1]^2$. The approximant is evaluated on a uniform
grid of $s = 120 \times 120 = 14400$ points. The number of PU subdomains
per direction is $n_{\mathrm{PU}} = \lfloor 65/2 \rfloor = 32$, giving
$K = n_{\mathrm{PU}}^2 = 1024$ patches with radius
$\delta_{\mathrm{PU}} = \sqrt{2}/n_{\mathrm{PU}}$.

\paragraph{Radial Basis Function.}
All methods use the Mat\'{e}rn $C^4$ kernel as RBF in \cref{rbf-ker},
which is mathematically expressed as:
\begin{align}
  \Phi(\varepsilon r) = \left(3 + 3\varepsilon r + \varepsilon^2 r^2\right)
  e^{-\varepsilon r},
\end{align}
with shape parameter $\varepsilon > 0$ controlling the flatness of the
kernel. For NL-PUM, $\varepsilon$ is fixed at $1$; for the remaining
methods, it is selected adaptively for each patch.

\paragraph{Weight Function.}
The compactly supported weight functions $\varphi_j$ appearing in
\cref{ShepWei,eq:nlweight} are taken as the Wendland $C^4$ RBFs:

\begin{align}
  \varphi_j(\bm{x}) &:= \varphi\!\left(1/\tilde{\rho}_j,
  \|\bm{x} - \bm{c}_j\|\right), \quad \text{where} \notag\\
  \varphi(\rho, r) &= \left(\max\!\left(1 - \rho\, r,\, 0\right)\right)^6
  \left(35(\rho\, r)^2 + 18(\rho\, r) + 3\right),
\end{align}

where $\tilde{\rho}_j$ is the (possibly shrunk) patch radius defined in
\cref{sec:radius_shrink}, with $\rho = 1/\tilde{\rho}_j$ its inverse. For
Linear PUM and LOOCV-NL-PUM, no shrinkage is applied and
$\tilde{\rho}_j = \delta_{\mathrm{PU}}$ for all $j$. These functions have
compact support within each patch and vanish smoothly at the patch
boundary, ensuring $C^4$ regularity of the PU.

\paragraph{LOOCV Parameters.}
For methods employing adaptive shape parameter selection, the search
interval is $\mathcal{I} = [0, \varepsilon_{\max}]$ with
$\varepsilon_{\max} = 10$. The LOOCV-NL-PUM-GOOI algorithm uses stopping
tolerance $\Delta = 0.03$, with $n_{\max}^{(1)} = n_{\max}^{(2)} = 10$
trial evaluations in the initial exploration and refinement stages,
respectively (see \cite{cav25b}). The smoothness threshold is set to
$\tau_\sigma = 0$, meaning that shrinkage is triggered for any patch with
a nonzero smoothness indicator; in practice the exponential form of
\cref{eq:rad_mod} ensures that the shrinkage is negligible when $\sigma_j$
is small (smooth regions) and significant only when $\sigma_j$ is large
(discontinuous regions). For LOOCV-NL-PUM, the patch radius is fixed at
$\delta_{\mathrm{PU}} = \sqrt{2}/n_{\mathrm{PU}}$ throughout. For
LOOCV-NL-PUM-GOOI, the initial radius follows \cref{eq:init_radius} with
expansion as described in \cref{sec:patch_init}, and the
discontinuity-aware shrinkage uses $\alpha = 3.5$ and non-linear exponent
$z = 4$. The minimum patch cardinality is $n_{\min} = 20$ for all methods.
Unless otherwise stated, all parameters listed in this section apply
uniformly to both the order-of-accuracy experiments of
\cref{sec:order_accuracy} and the discontinuous test cases of
\cref{sec:results_discussion}.

\subsection{Error Metrics}

The error measures below are consistent with those of
\cref{sec:order_accuracy}, but with the evaluation set $\hat{X}$
corresponding to the uniform grid consisting of $s = 120\times120 = 14400$
points as expressed in \cref{setup}. We report two error measures. The
Maximum Absolute Error (MAE) is defined as:
\begin{align}
  \text{MAE} = \max_{1 \leq k \leq s}
  \bigl|f(\hat{\bm{x}}_k) - I(v)(\hat{\bm{x}}_k)\bigr|,
\end{align}
and the discrete $\ell^2$ error norm is:
\begin{align}
  \ell^2\;\text{error} = \sqrt{\frac{1}{s} \sum_{k=1}^{s}
  \bigl|f(\hat{\bm{x}}_k) - I(v)(\hat{\bm{x}}_k)\bigr|^2},
\end{align}
where $\hat{\bm{x}}_k \in \hat{X}$ are the evaluation points defined in
\cref{setup}, and $I(v)$ denotes the PU-based approximant of the method
under consideration. CPU times are measured using MATLAB \texttt{tic}/\texttt{toc}
commands and reported in seconds.

\subsection{Results}
\label{sec:results_discussion}

\begin{figure}[htbp]
\centering
\includegraphics[width=0.48\textwidth]{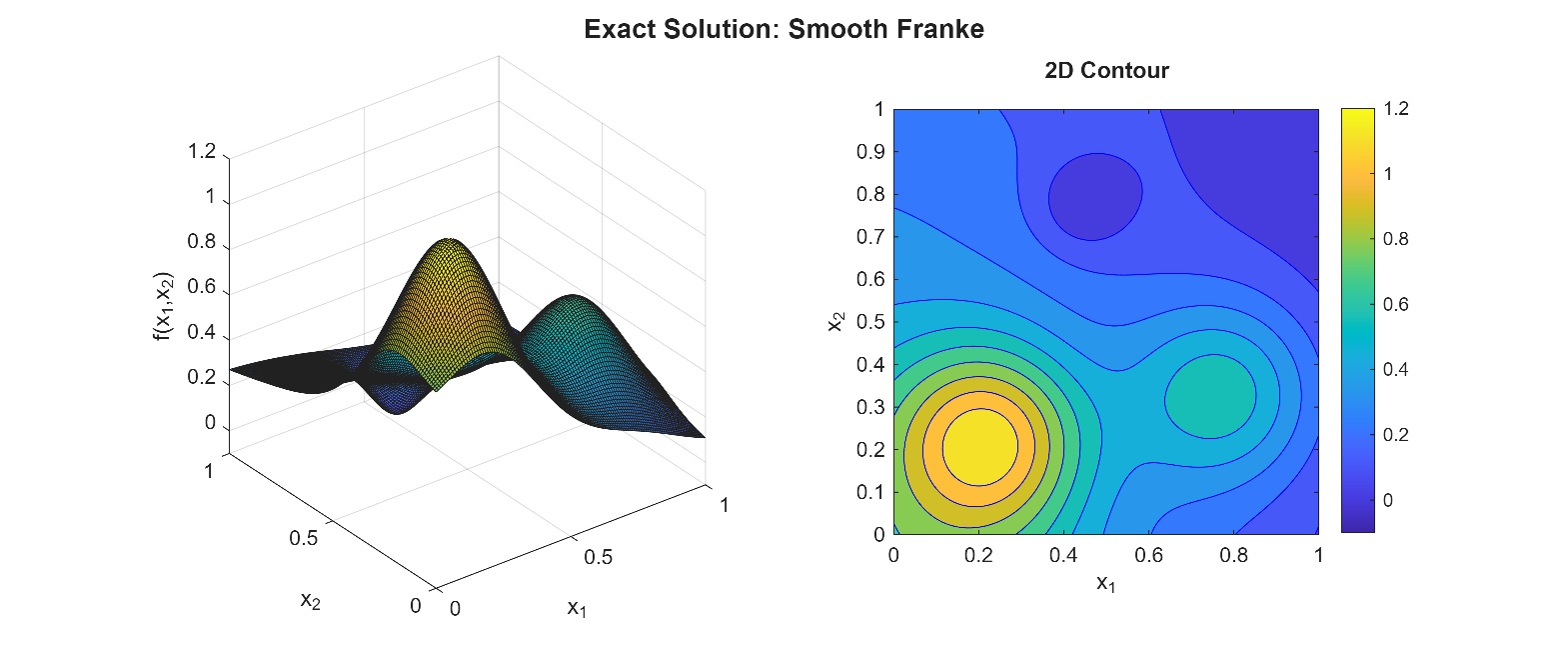}\hfill
\includegraphics[width=0.48\textwidth]{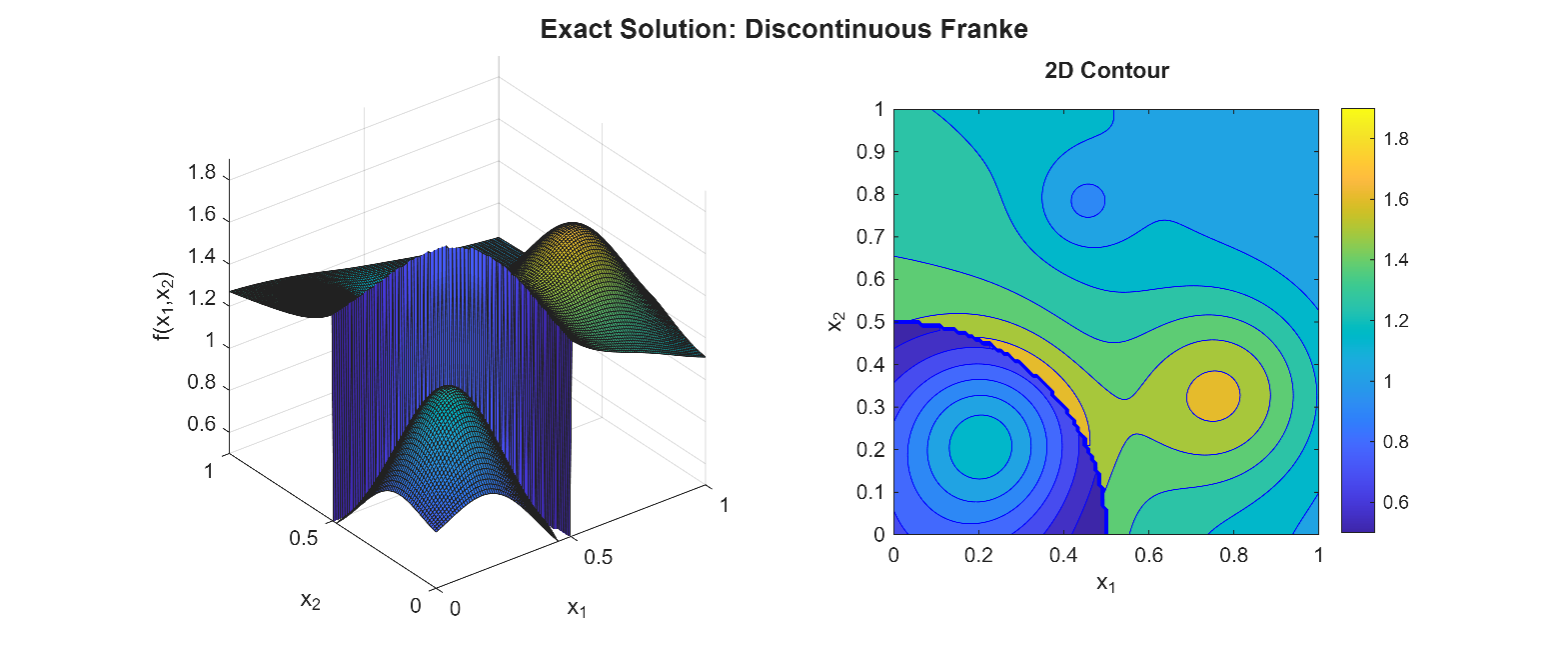}

\vspace{0.5em}
\includegraphics[width=0.48\textwidth]{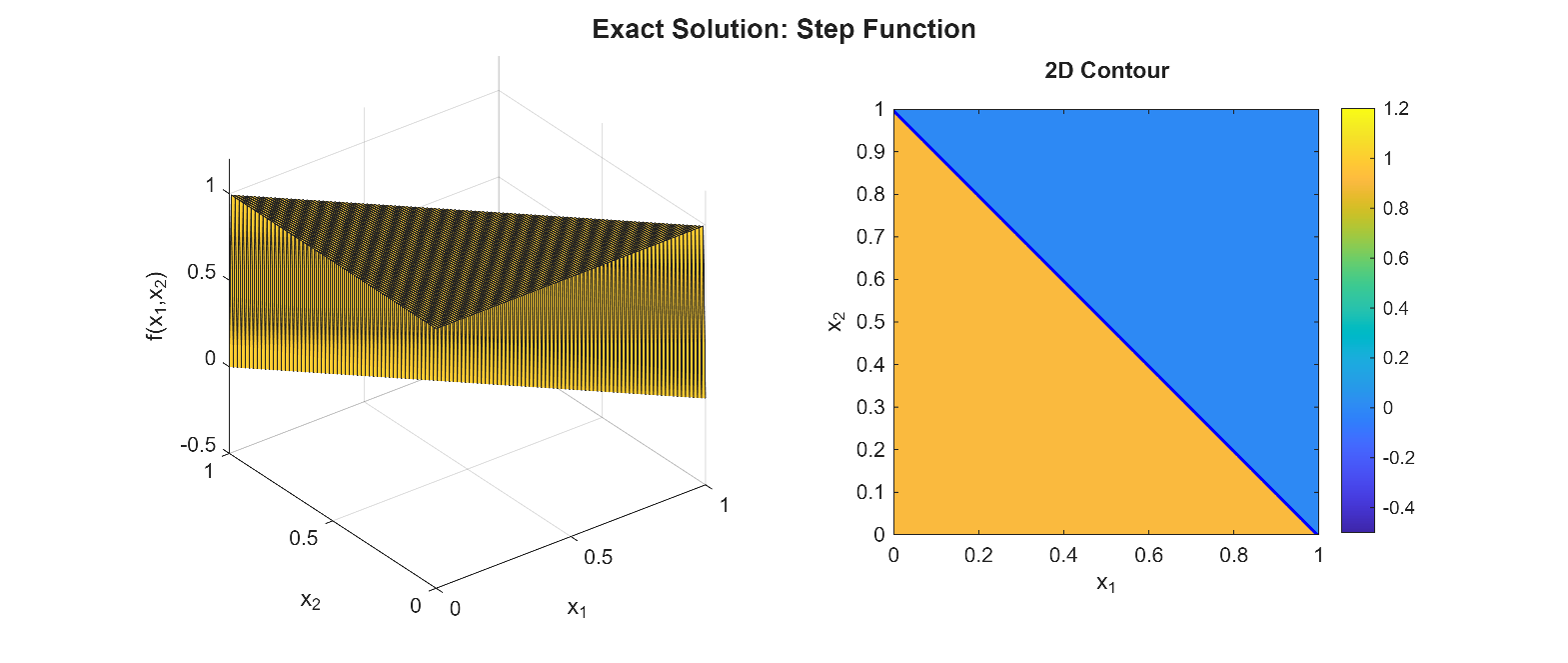}\hfill
\includegraphics[width=0.48\textwidth]{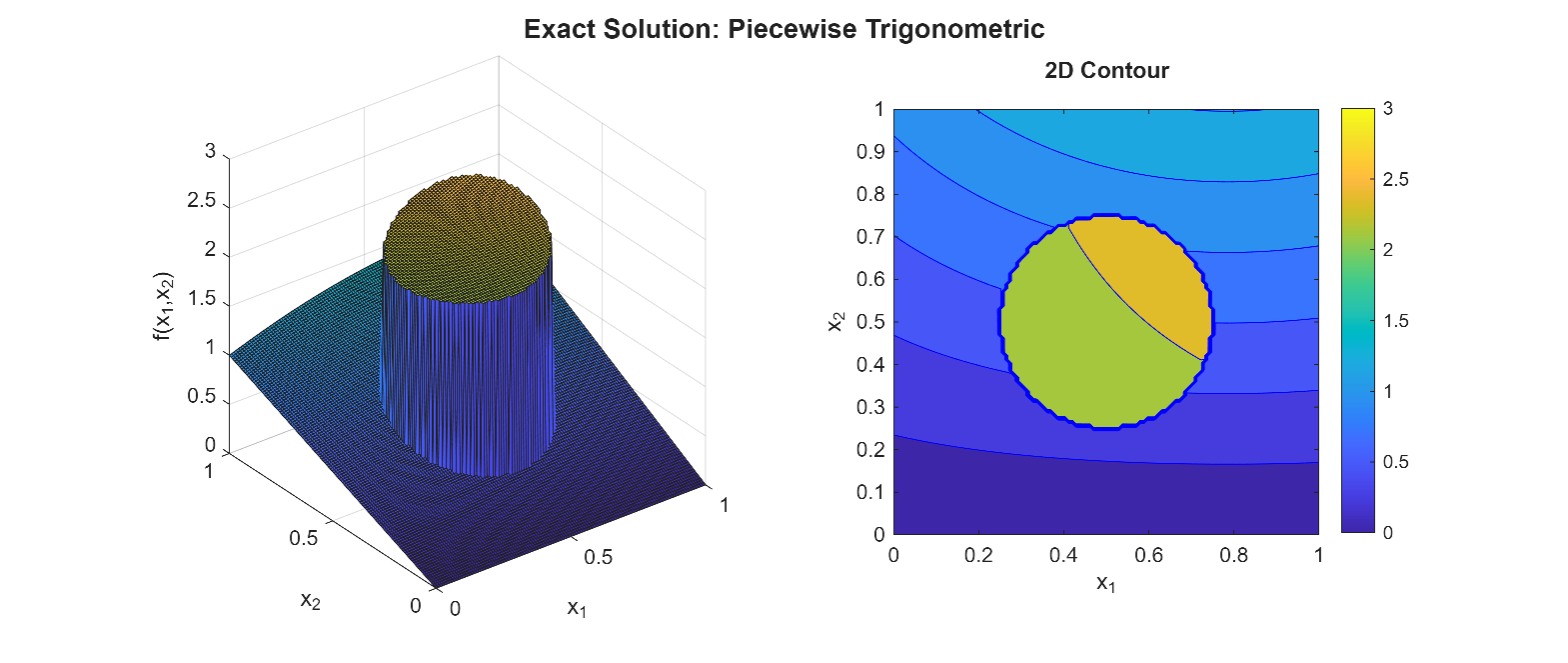}
\caption{Exact solutions (3D surface and 2D contour) for the four test
functions. Top left: smooth Franke function $f_1$. Top right:
discontinuous Franke function $f_2$. Bottom left: step function $f_3$.
Bottom right: piecewise trigonometric function $f_4$.}
\label{fig:all_exact_solutions}
\end{figure}

\Cref{tab:pu_results} reports MAE, $\ell^2$ error, and CPU time for all
four methods. The exact solutions are shown in
\cref{fig:all_exact_solutions}, while comparisons among the methods are
provided in \cref{fig:methods_smooth_franke,fig:methods_disc_franke,%
fig:methods_step,fig:methods_trig}.

\begin{figure}[htbp]
\centering
\includegraphics[width=\textwidth]{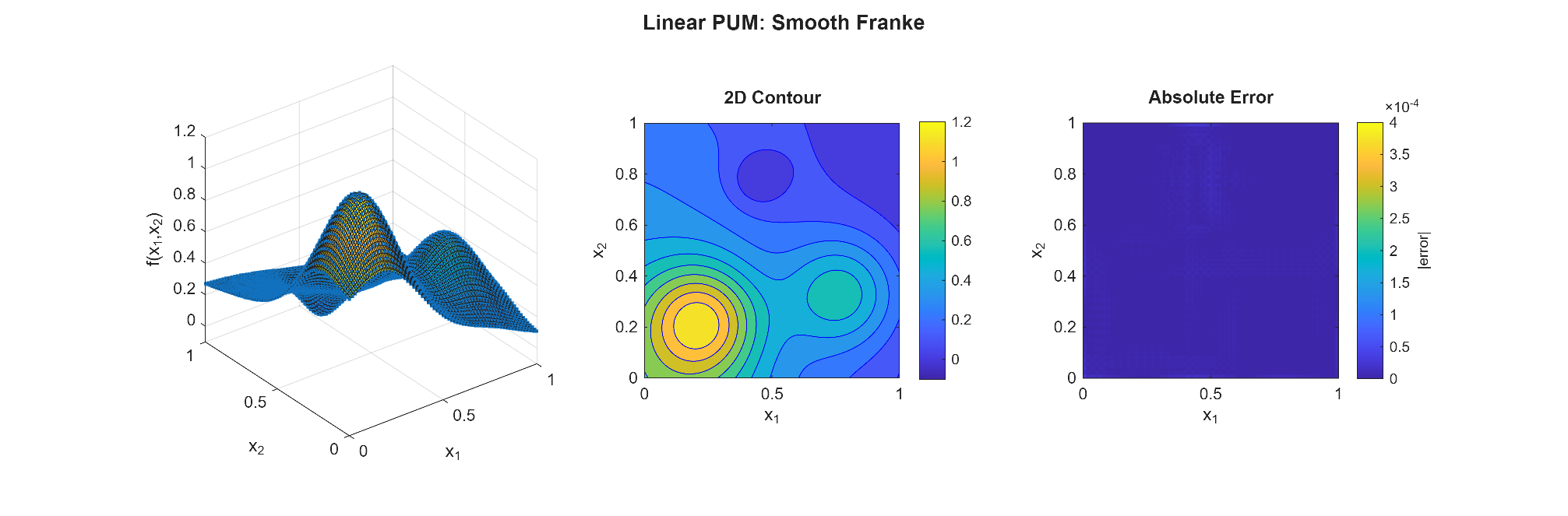}

\vspace{0.5em}
\includegraphics[width=\textwidth]{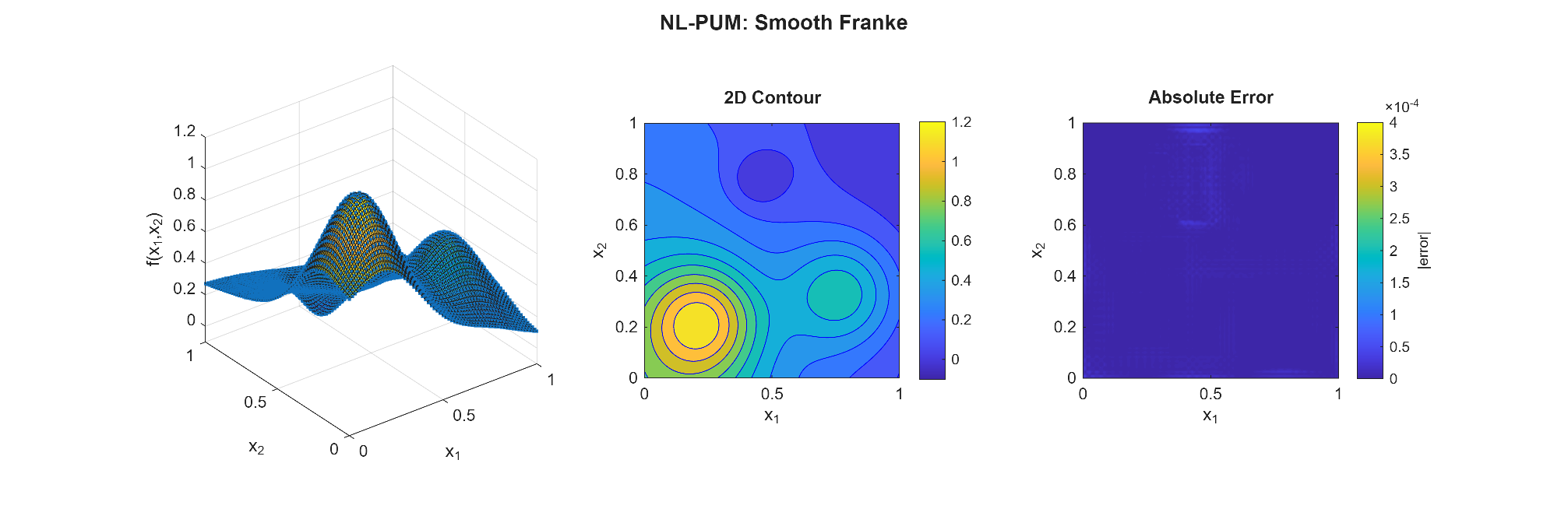}

\vspace{0.5em}
\includegraphics[width=\textwidth]{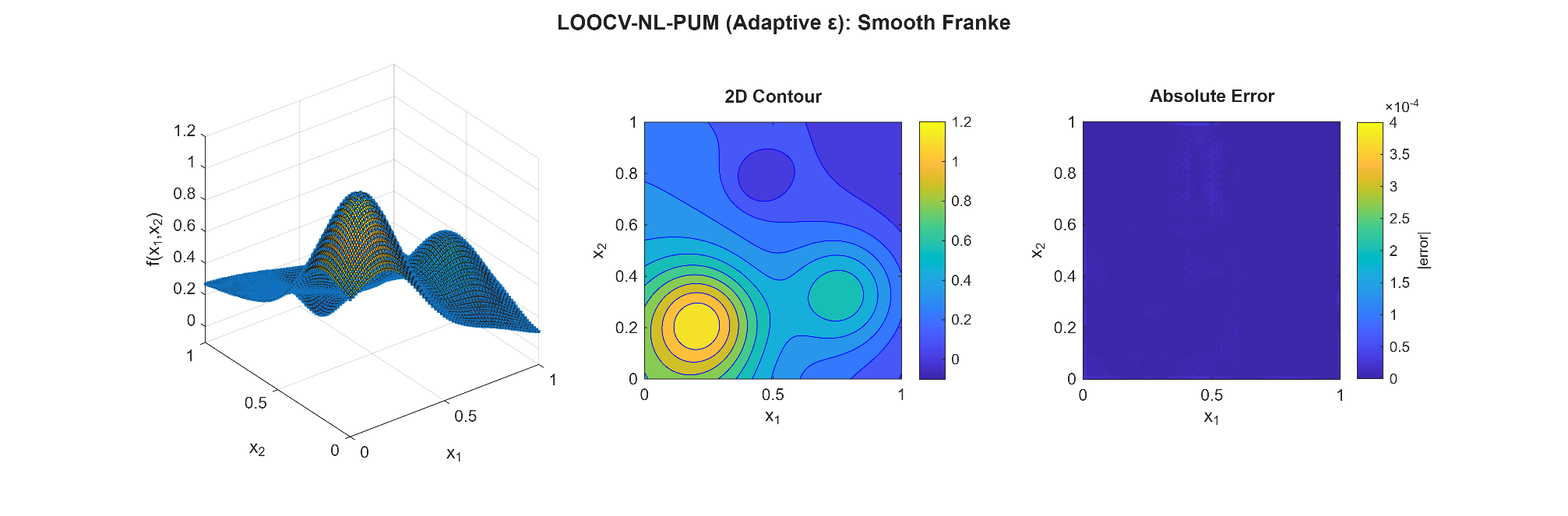}

\vspace{0.5em}
\includegraphics[width=\textwidth]{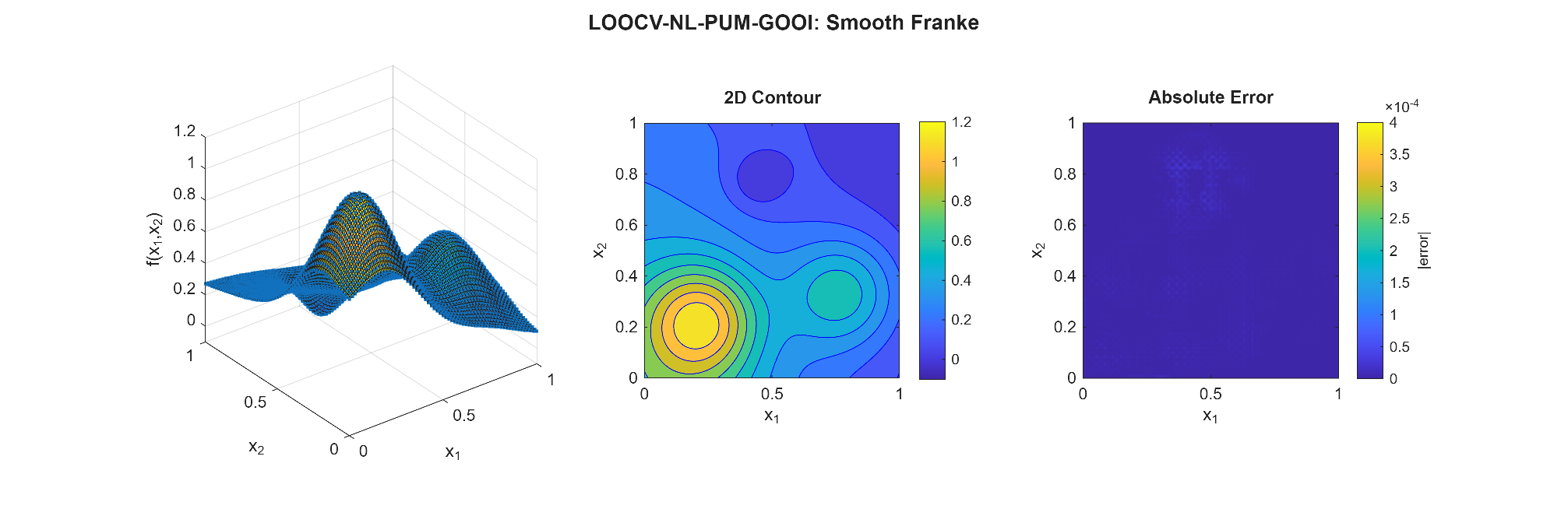}
\caption{Comparison of all methods for the smooth Franke function $f_1$.}
\label{fig:methods_smooth_franke}
\end{figure}

In the smooth Franke function $f_1$, all methods perform comparably, as
\cref{fig:methods_smooth_franke} confirms visually. LOOCV-NL-PUM-GOOI
achieves the best errors (MAE $= 3.38 \times 10^{-5}$,
$\ell^2 = 3.95 \times 10^{-6}$), but the differences are modest and no
shrinkage is triggered, so this case serves primarily as a sanity check
that the proposed method does not introduce an accuracy penalty in smooth
regions.

\begin{figure}[htbp]
\centering
\includegraphics[width=\textwidth]{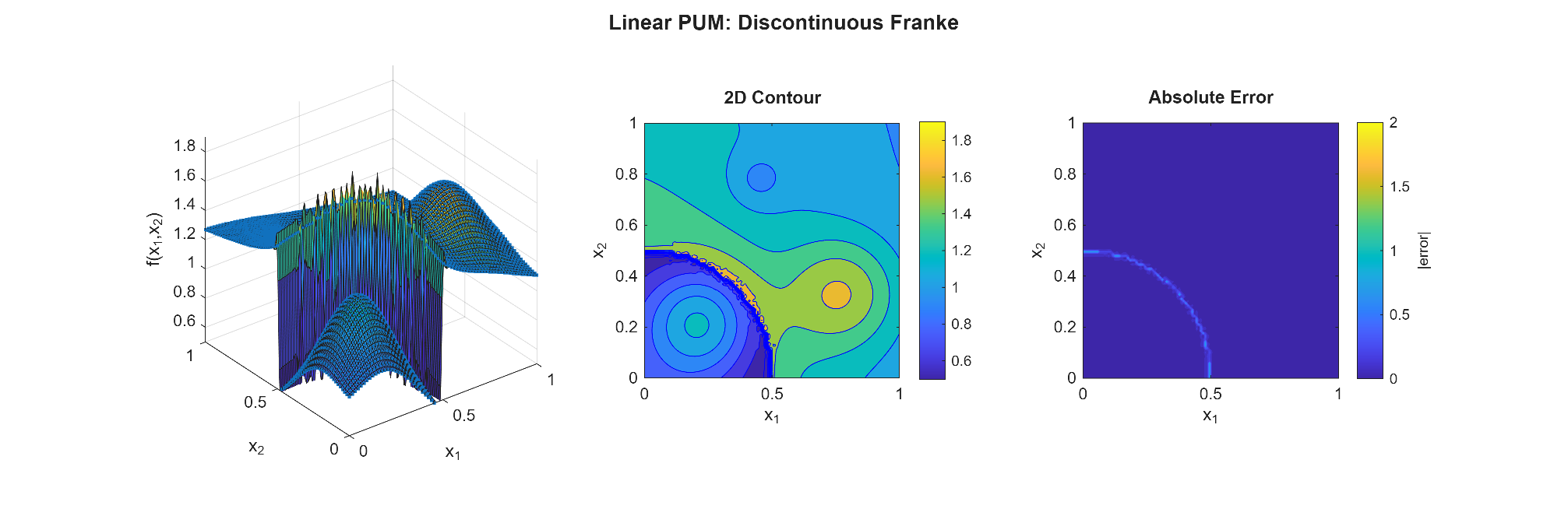}

\vspace{0.5em}
\includegraphics[width=\textwidth]{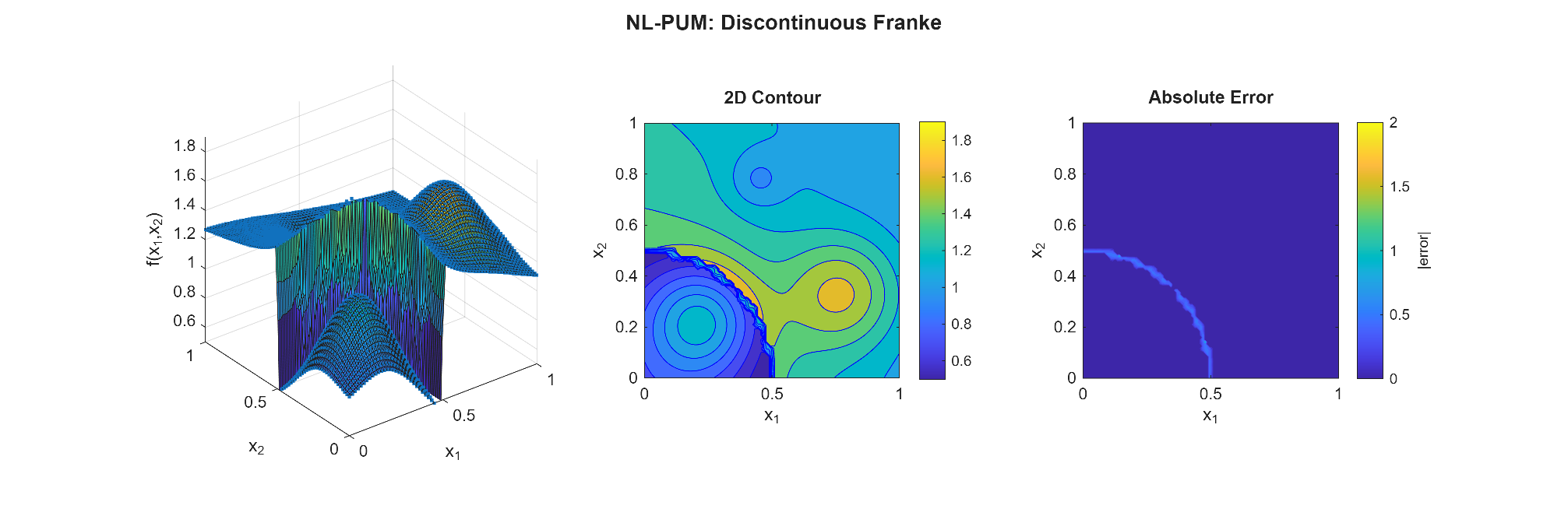}

\vspace{0.5em}
\includegraphics[width=\textwidth]{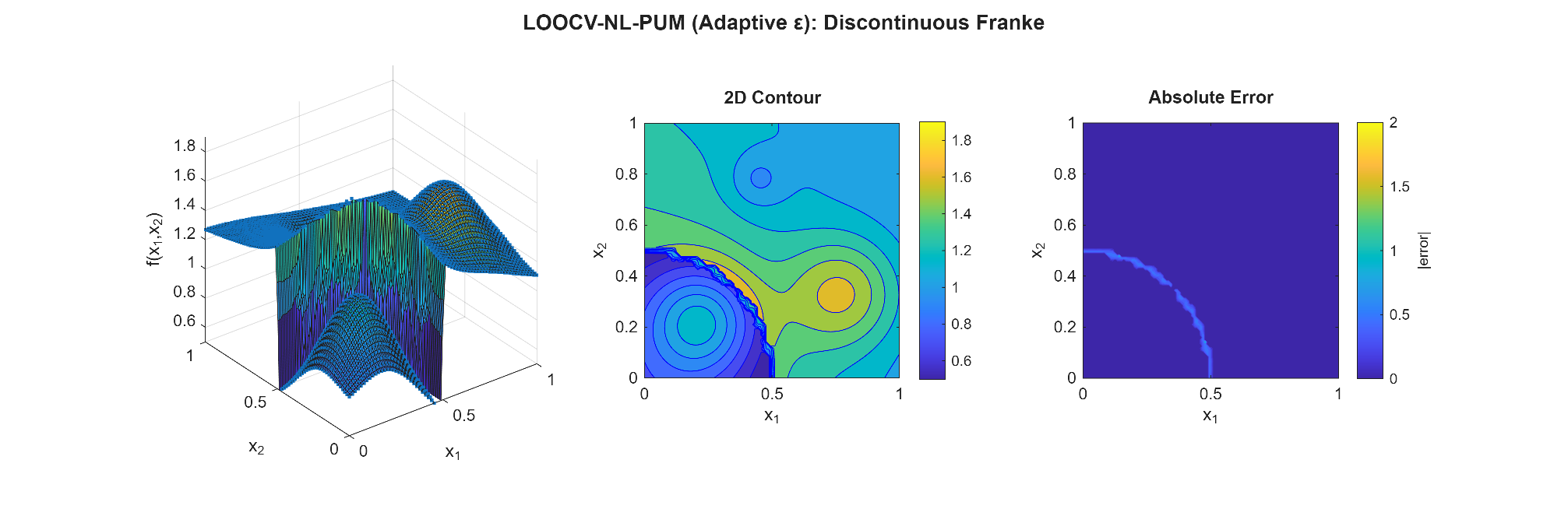}

\vspace{0.5em}
\includegraphics[width=\textwidth]{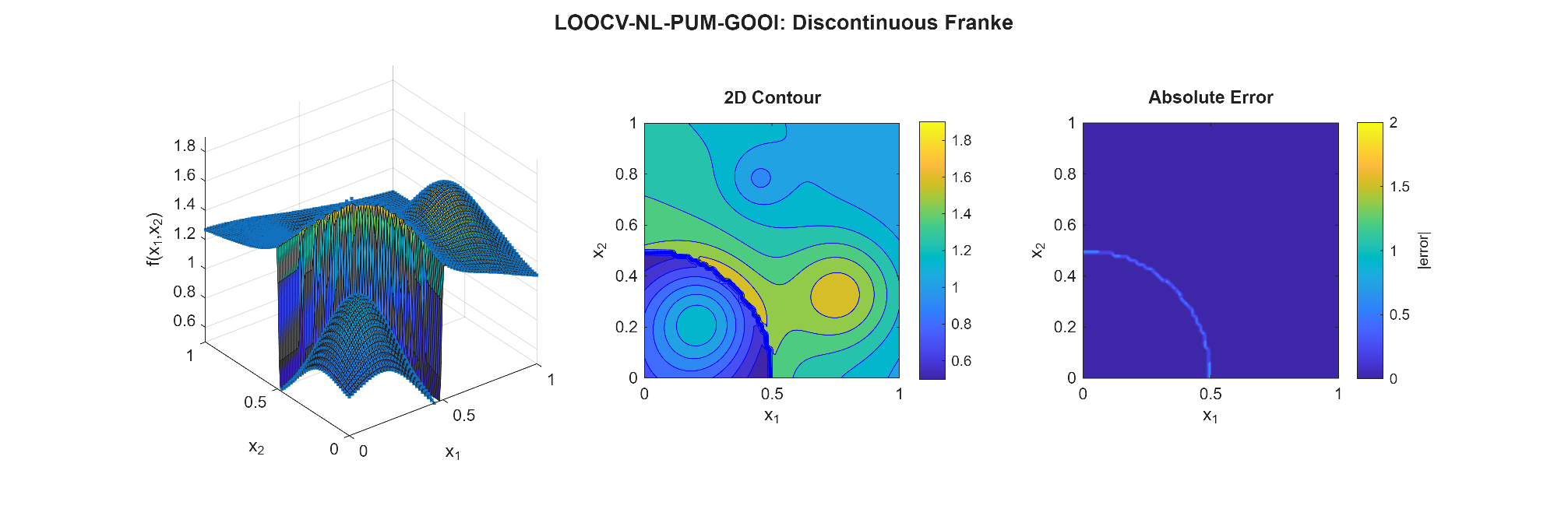}
\caption{Comparison of all methods for the discontinuous Franke function
$f_2$.}
\label{fig:methods_disc_franke}
\end{figure}

\begin{figure}[htbp]
\centering
\includegraphics[width=\textwidth]{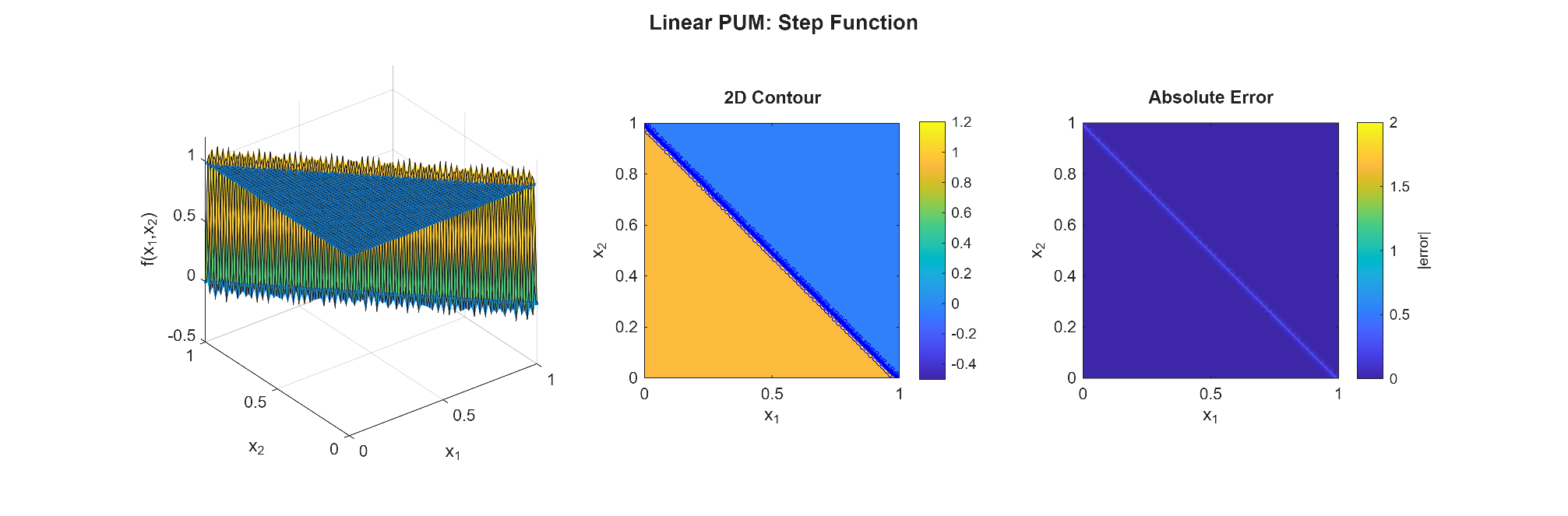}

\vspace{0.5em}
\includegraphics[width=\textwidth]{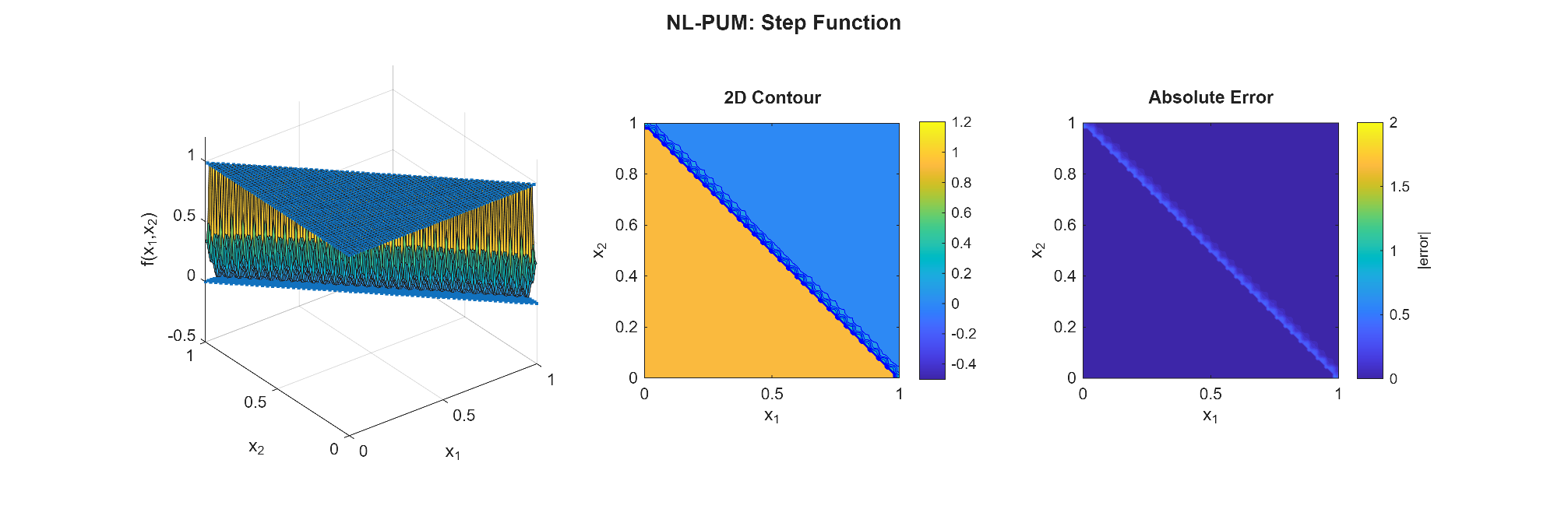}

\vspace{0.5em}
\includegraphics[width=\textwidth]{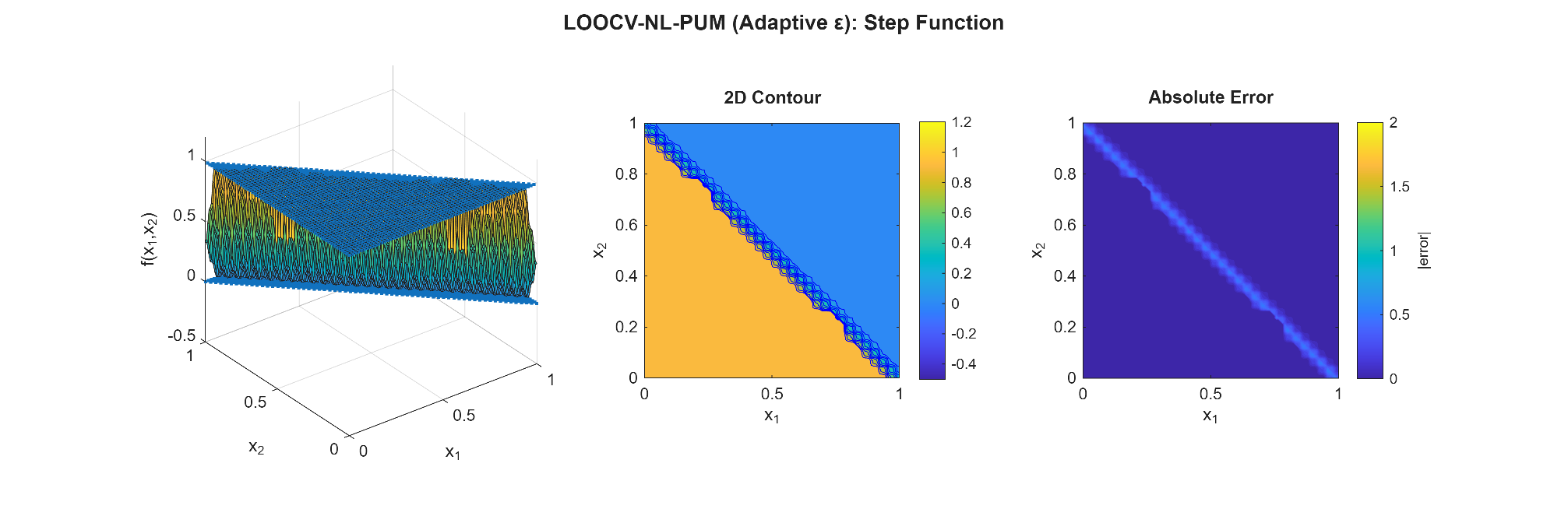}

\vspace{0.5em}
\includegraphics[width=\textwidth]{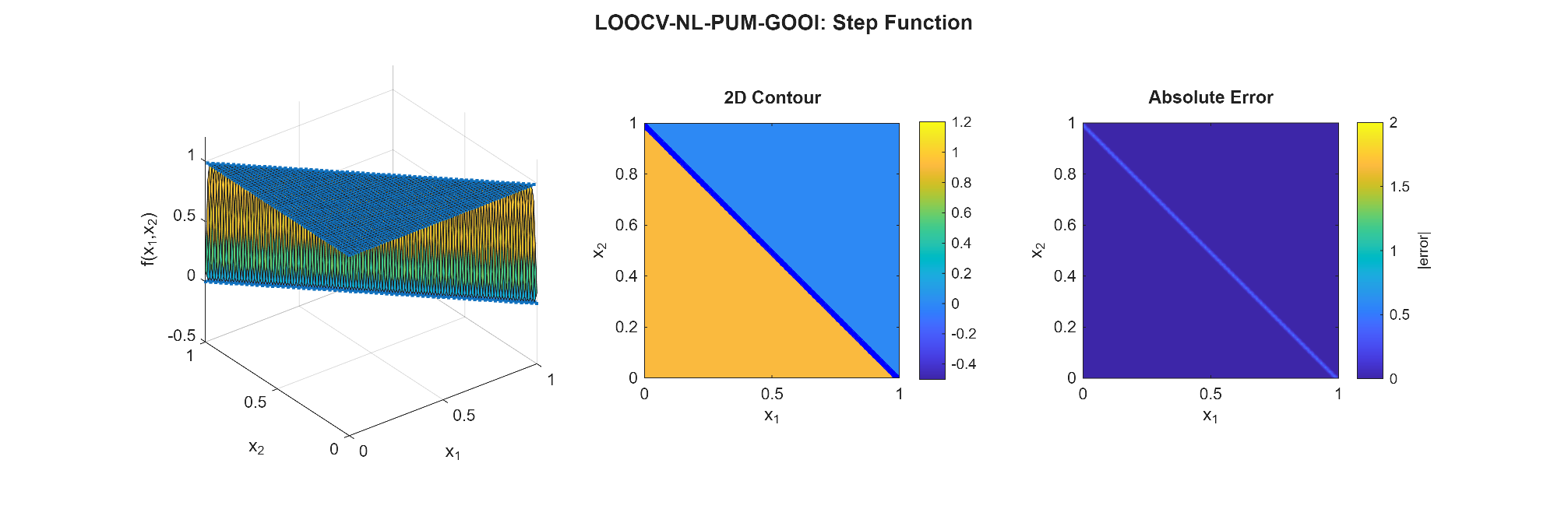}
\caption{Comparison of all methods for the step function $f_3$.}
\label{fig:methods_step}
\end{figure}

\begin{figure}[htbp]
\centering
\includegraphics[width=\textwidth]{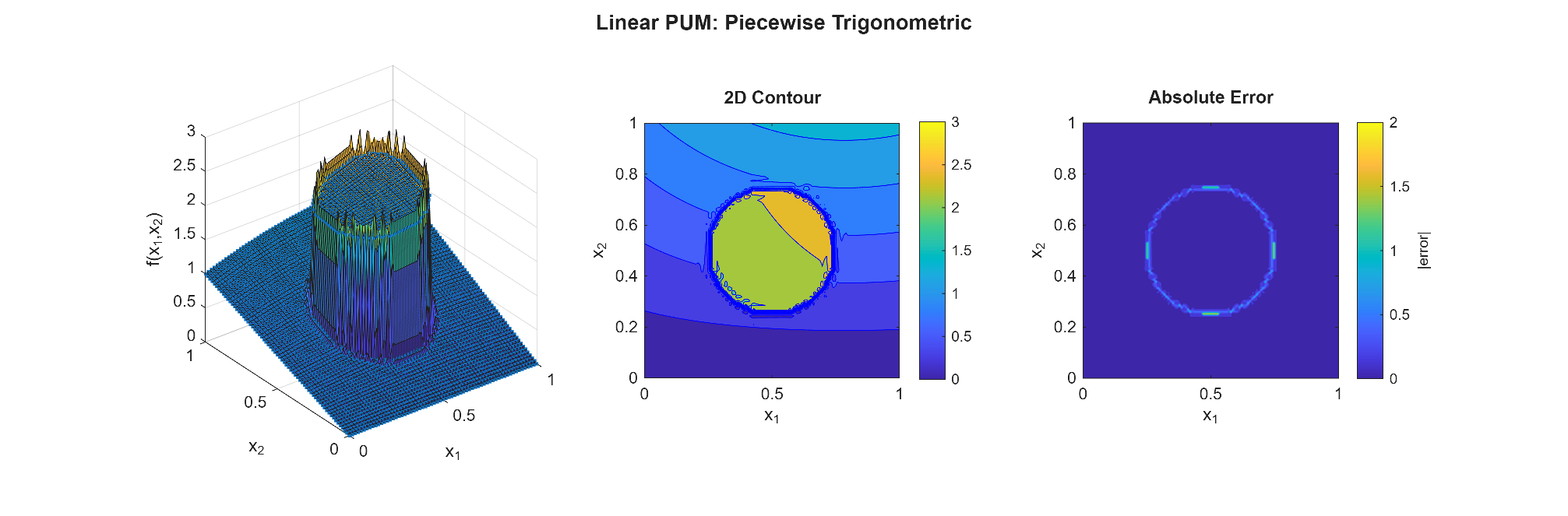}

\vspace{0.5em}
\includegraphics[width=\textwidth]{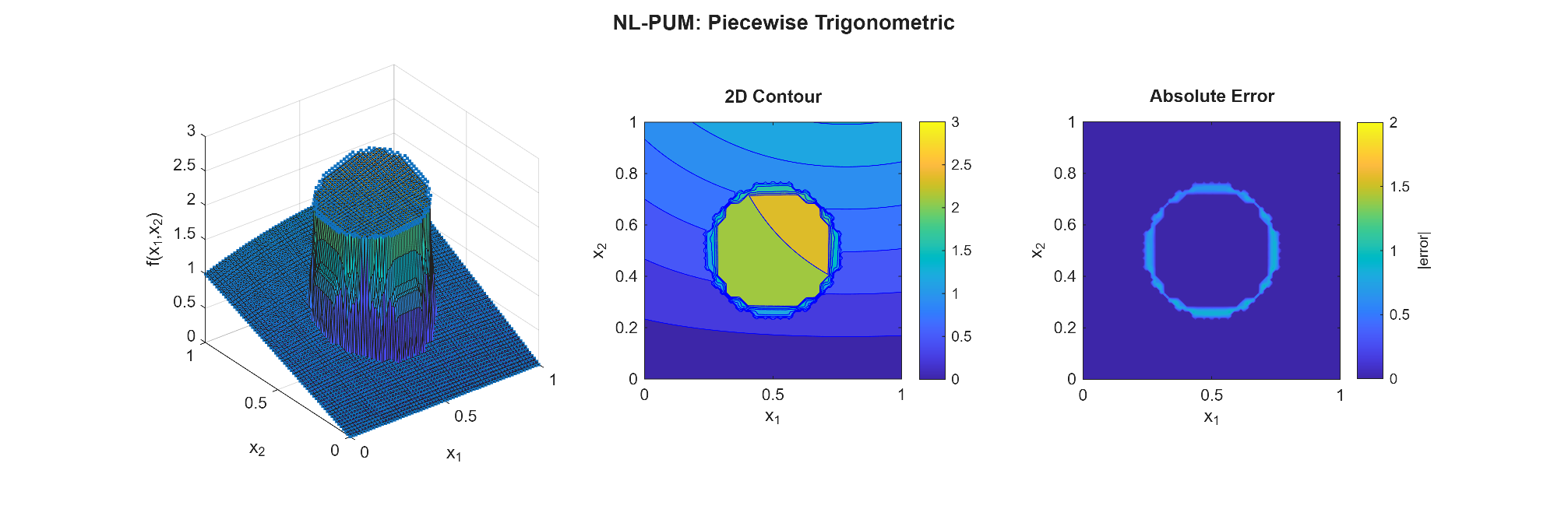}

\vspace{0.5em}
\includegraphics[width=\textwidth]{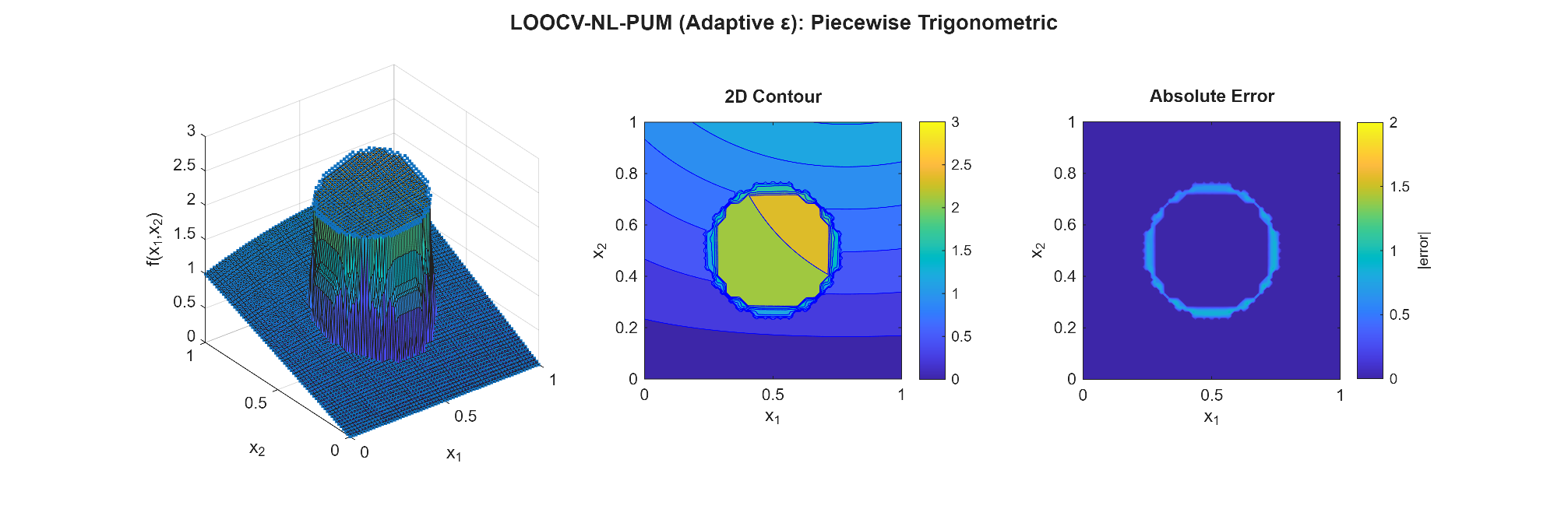}

\vspace{0.5em}
\includegraphics[width=\textwidth]{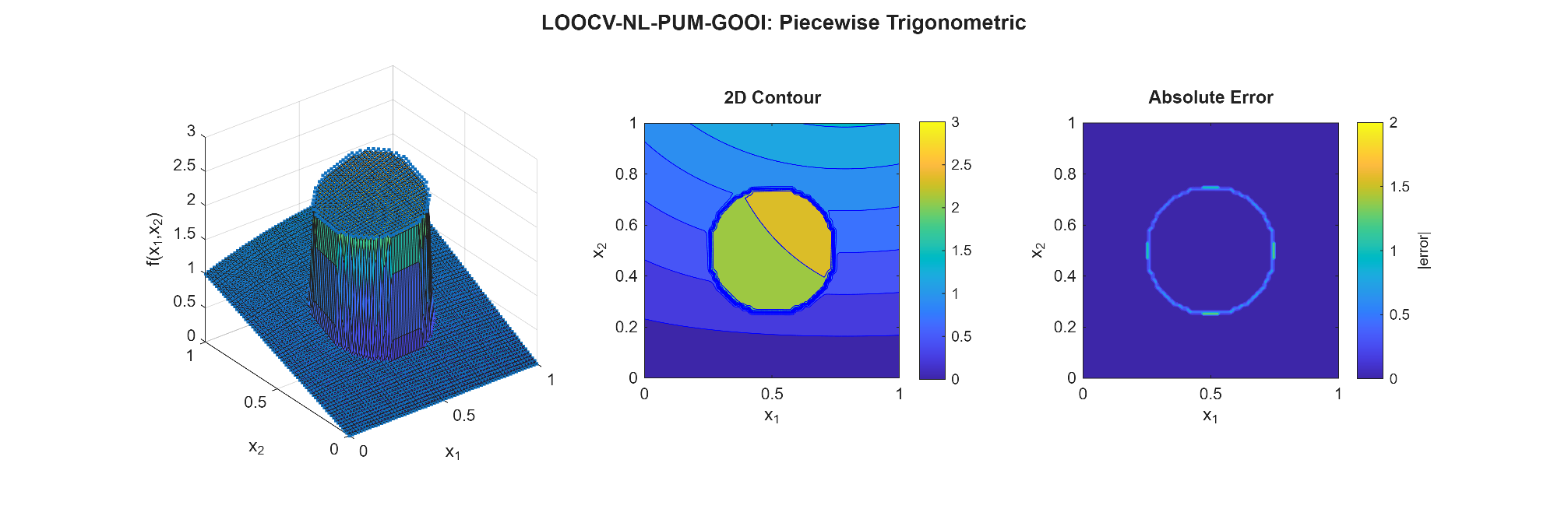}
\caption{Comparison of all methods for the piecewise trigonometric function
$f_4$.}
\label{fig:methods_trig}
\end{figure}

The discontinuous test cases expose a much clearer picture. Across all
three---the discontinuous Franke function $f_2$
(\cref{fig:methods_disc_franke}), the step function $f_3$
(\cref{fig:methods_step}), and the piecewise trigonometric function $f_4$
(\cref{fig:methods_trig})---the most striking observation from the absolute
error panels is that LOOCV-NL-PUM is indistinguishable from NL-PUM:
adapting the shape parameter $\varepsilon$ alone, while keeping the patch
radius fixed, produces no perceptible reduction in the oscillations near
the interface. This confirms that the radius itself, not the kernel shape,
is the critical quantity governing behaviour near discontinuities. Linear
PUM produces the most severe oscillations in every case.
LOOCV-NL-PUM-GOOI, on the other hand, yields a visibly cleaner surface in
each figure: the oscillatory curtains along the circular interface
$\mathcal{C}$ and the diagonal interface $\mathcal{L}$ are substantially
reduced, and the flat regions on either side of the jump are recovered more
faithfully. These visual improvements are consistent with the quantitative
results in \cref{tab:pu_results}, where LOOCV-NL-PUM-GOOI achieves the
best $\ell^2$ error in every discontinuous test case. The one exception is
the MAE on the step and piecewise trigonometric functions, where radius
shrinkage near a compact interior interface can leave isolated evaluation
points with insufficient patch coverage, elevating the maximum pointwise
error at those locations while the global $\ell^2$ quality improves. This
is a known trade-off of radius-shrinkage strategies and motivates future
geometry-aware refinements of the shrinkage criterion.

Regarding computational cost, LOOCV-NL-PUM-GOOI ($3.47$--$3.58$~s) is
comparable to NL-PUM ($2.79$--$3.12$~s) and faster than LOOCV-NL-PUM
($5.31$--$5.64$~s), since radius shrinkage reduces local dataset sizes near
discontinuities and partially offsets the overhead of the adaptive
parameter search. Overall, LOOCV-NL-PUM-GOOI provides the best
accuracy-to-cost ratio among the adaptive methods tested.

\begin{table}[htbp]
\caption{Performance comparison of PU methods across all test functions.}
\label{tab:pu_results}
\begin{center}
\resizebox{\textwidth}{!}{%
\begin{tabular}{llccc}
\toprule
\textbf{Test Function} & \textbf{Method} & \textbf{MAE}
& $\bm{\ell^2}$ \textbf{Error} & \textbf{CPU (s)} \\
\midrule
\multirow{4}{*}{\shortstack[c]{Smooth Franke\\$f_1$}}
& Linear PUM                 & 1.70e$-$04 & 5.60e$-$06 & 5.13e$-$01 \\
& NL-PUM                     & 3.74e$-$05 & 6.56e$-$06 & 2.79e$+$00 \\
& LOOCV-NL-PUM               & 3.61e$-$05 & 6.89e$-$06 & 5.64e$+$00 \\
& \textbf{LOOCV-NL-PUM-GOOI}
  & \textbf{3.38e$-$05} & \textbf{3.95e$-$06} & 3.53e$+$00 \\
\midrule
\multirow{4}{*}{\shortstack[c]{Discontinuous\\Franke $f_2$}}
& Linear PUM                 & 1.19e$+$00 & 5.35e$-$02 & 4.60e$-$01 \\
& NL-PUM                     & 1.00e$+$00 & 1.01e$-$01 & 3.12e$+$00 \\
& LOOCV-NL-PUM               & 1.00e$+$00 & 1.01e$-$01 & 5.56e$+$00 \\
& \textbf{LOOCV-NL-PUM-GOOI}
  & \textbf{7.42e$-$01} & \textbf{4.41e$-$02} & 3.47e$+$00 \\
\midrule
\multirow{4}{*}{\shortstack[c]{Step Function\\$f_3$}}
& Linear PUM                 & 1.38e$+$00 & 6.73e$-$02 & 4.77e$-$01 \\
& NL-PUM                     & 6.39e$-$01 & 1.15e$-$01 & 2.95e$+$00 \\
& LOOCV-NL-PUM               & 6.39e$-$01 & 1.72e$-$01 & 5.31e$+$00 \\
& \textbf{LOOCV-NL-PUM-GOOI}
  & \textbf{5.53e$-$01} & \textbf{4.87e$-$02} & 3.58e$+$00 \\
\midrule
\multirow{4}{*}{\shortstack[c]{Piecewise\\Trigonometric $f_4$}}
& Linear PUM                 & 1.91e$+$00 & 1.18e$-$01 & 4.68e$-$01 \\
& NL-PUM                     & 1.04e$+$00 & 2.67e$-$01 & 2.81e$+$00 \\
& LOOCV-NL-PUM               & 1.04e$+$00 & 2.67e$-$01 & 5.49e$+$00 \\
& \textbf{LOOCV-NL-PUM-GOOI}
  & 1.52e$+$00 & \textbf{1.02e$-$01} & 3.54e$+$00 \\
\bottomrule
\end{tabular}
}
\end{center}
\end{table}

We further validate the proposed method on the ETOPO~2022 global relief
model \cite{etopo2022}, using elevation data from the Norwegian Fjords
(Sognefjord region), covering the geographic window
$[\,6^{\circ}\text{E},\, 9^{\circ}\text{E}\,] \times
[\,60^{\circ}\text{N},\, 62^{\circ}\text{N}\,]$. This region presents
mountain peaks rising above $2800$~m alongside fjord channels descending
below $-1000$~m, producing multiple sharp elevation discontinuities that
are really hard to reconstruct using standard interpolation methods.

The full dataset consists of $21{,}600$ data points but to develop a
realistic reconstruction scenario we took only $3{,}600$ scattered points
that are less than $17\%$ of the whole dataset. For a fair comparison we
used the Mat\'{e}rn $C^4$ kernel and Wendland $C^4$ weights, with
$n_{\mathrm{PU}} = 30$ PU centers per direction,
$K = n_{\mathrm{PU}}^2 = 900$ patches, and a $s = 120 \times 120$
evaluation grid for all the methods. Since the true surface is unknown, no
quantitative error metrics are reported; results are visually assessed.

\begin{figure}[htbp]
\centering
\begin{minipage}{0.48\textwidth}
  \centering
  \includegraphics[width=\textwidth]{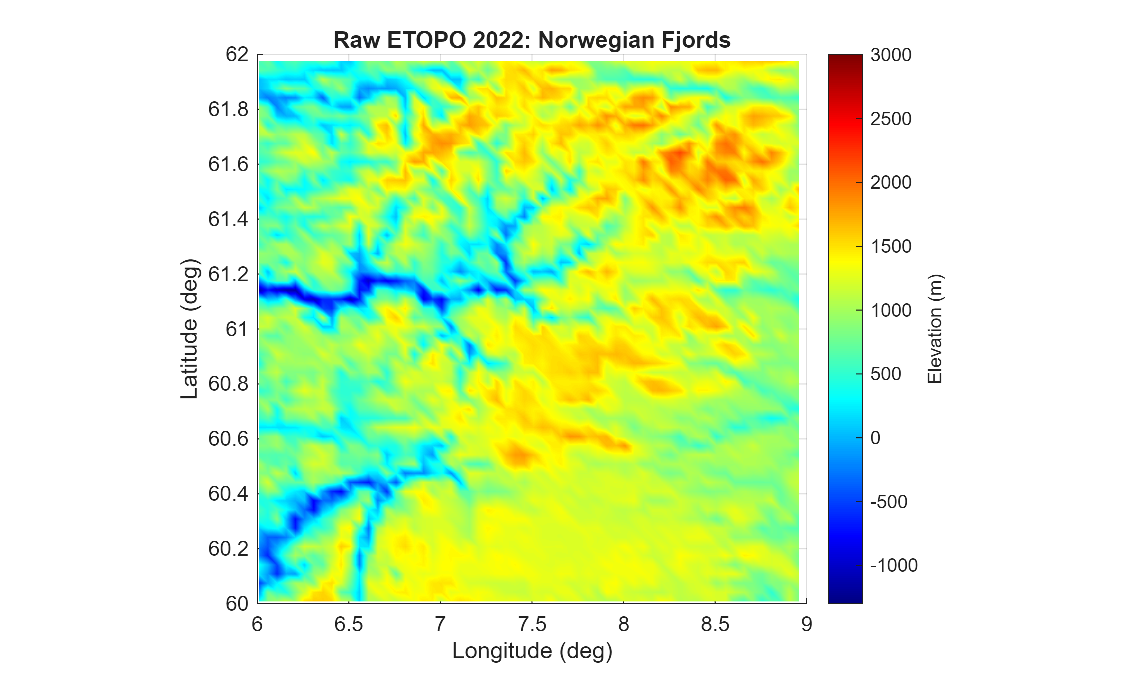}
\end{minipage}\hfill
\begin{minipage}{0.48\textwidth}
  \centering
  \includegraphics[width=\textwidth]{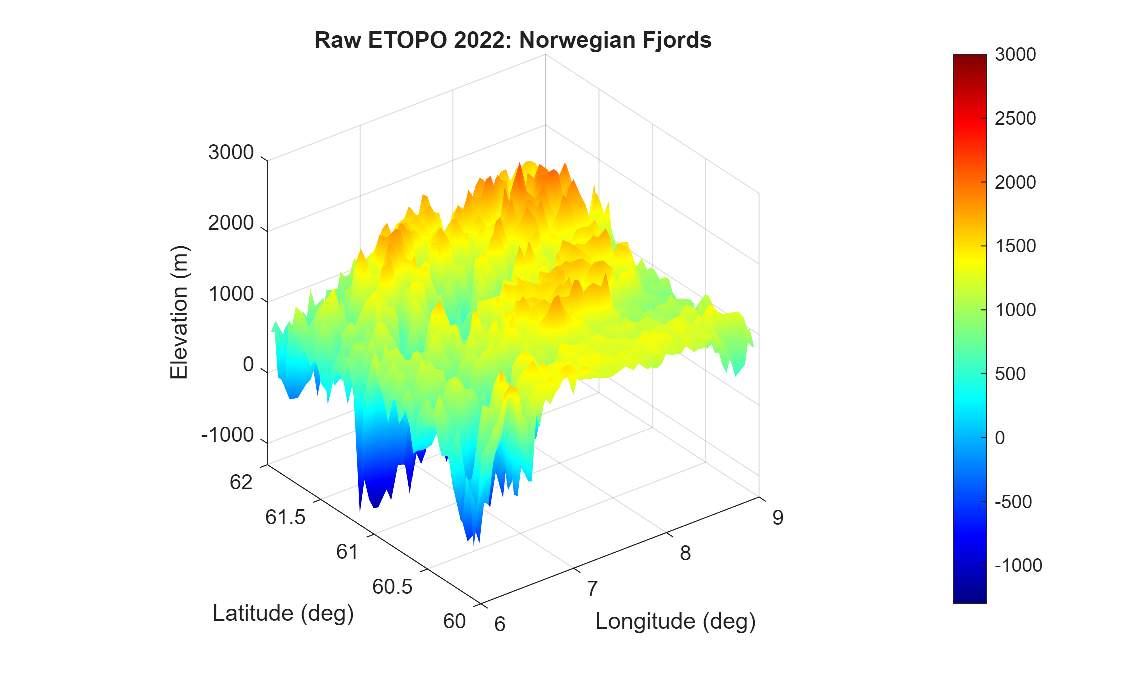}
\end{minipage}

\vspace{0.5em}
\begin{minipage}{0.48\textwidth}
  \centering
  \includegraphics[width=\textwidth]{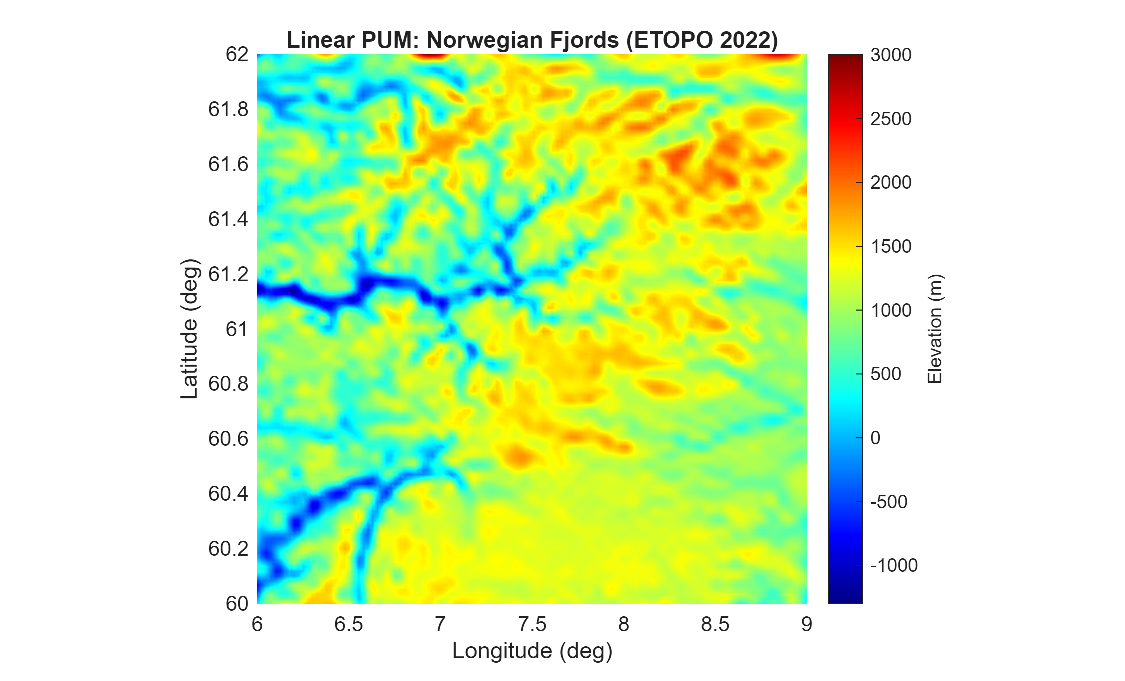}
\end{minipage}\hfill
\begin{minipage}{0.48\textwidth}
  \centering
  \includegraphics[width=\textwidth]{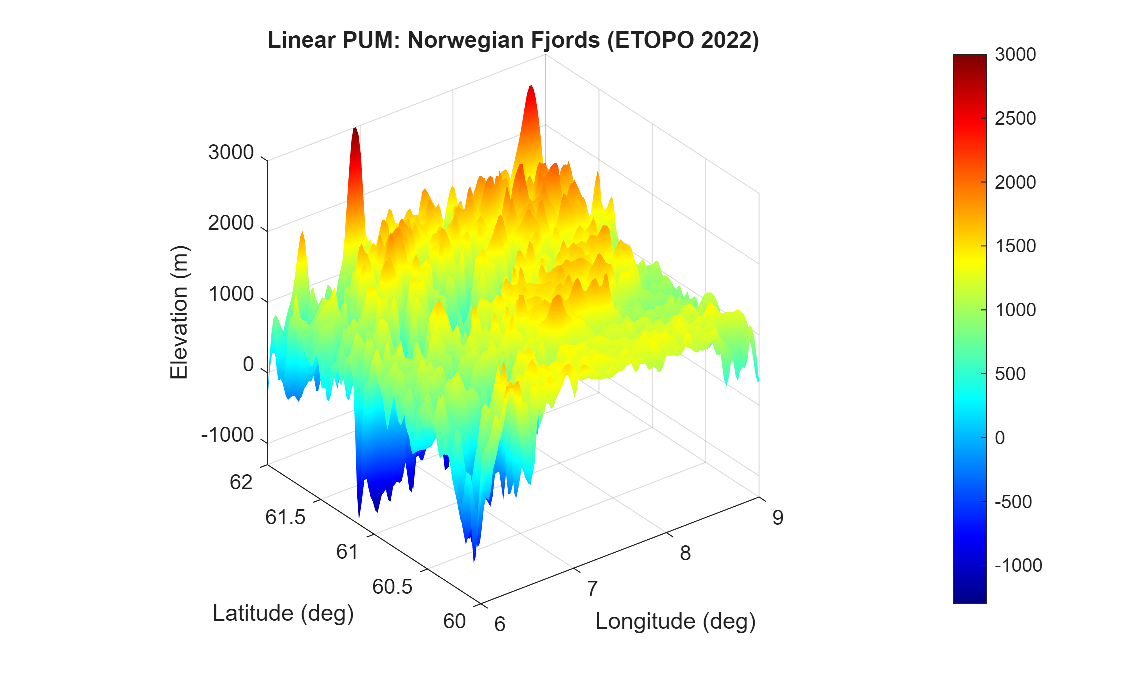}
\end{minipage}

\vspace{0.5em}
\begin{minipage}{0.48\textwidth}
  \centering
  \includegraphics[width=\textwidth]{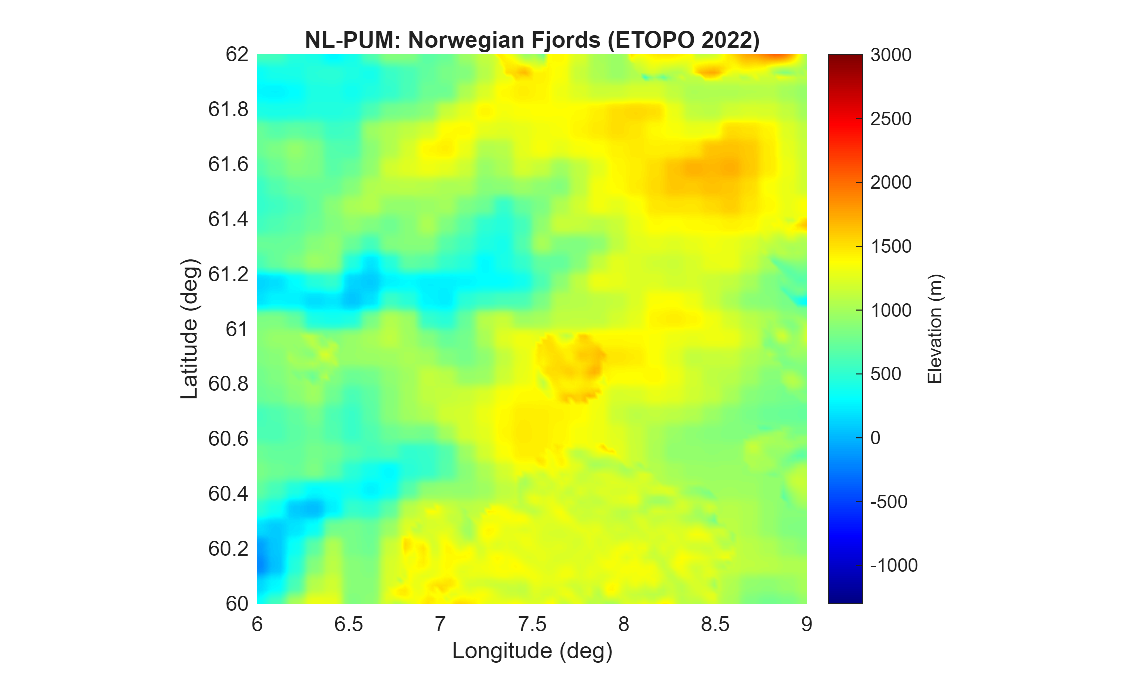}
\end{minipage}\hfill
\begin{minipage}{0.48\textwidth}
  \centering
  \includegraphics[width=\textwidth]{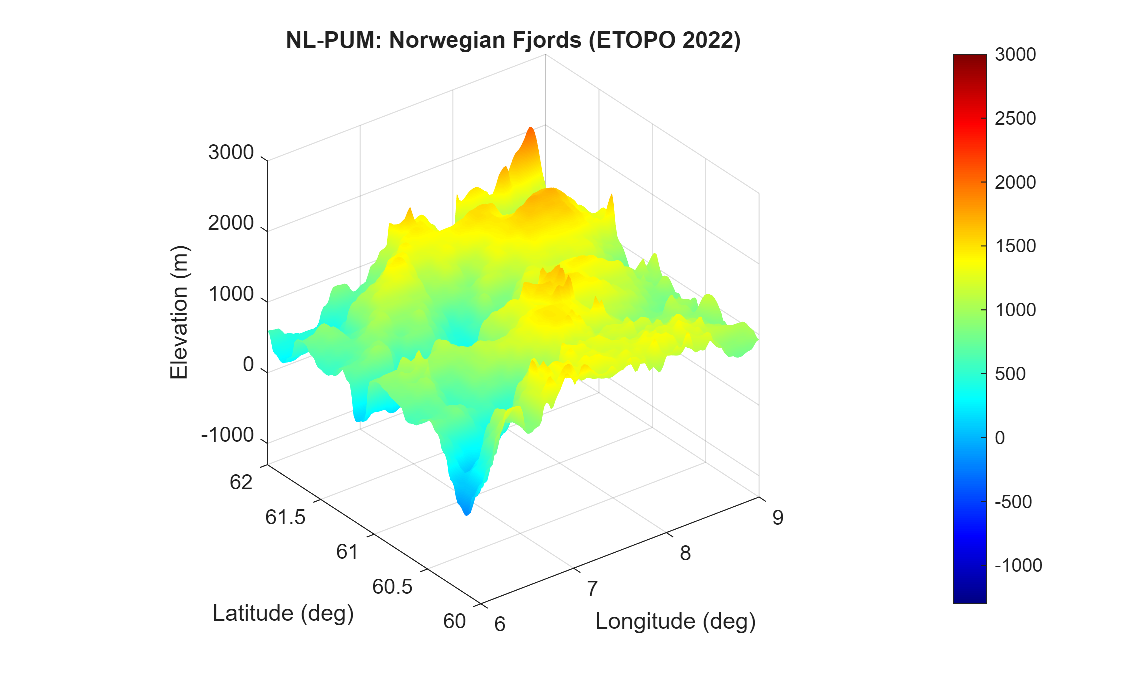}
\end{minipage}

\vspace{0.5em}
\begin{minipage}{0.48\textwidth}
  \centering
  \includegraphics[width=\textwidth]{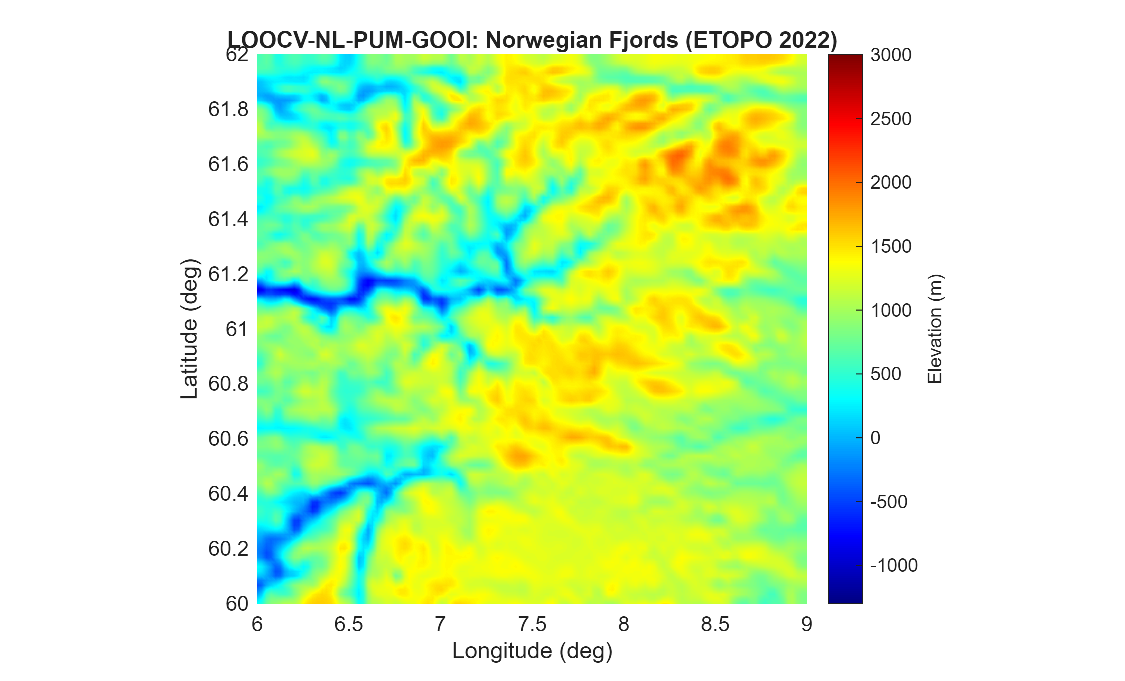}
\end{minipage}\hfill
\begin{minipage}{0.48\textwidth}
  \centering
  \includegraphics[width=\textwidth]{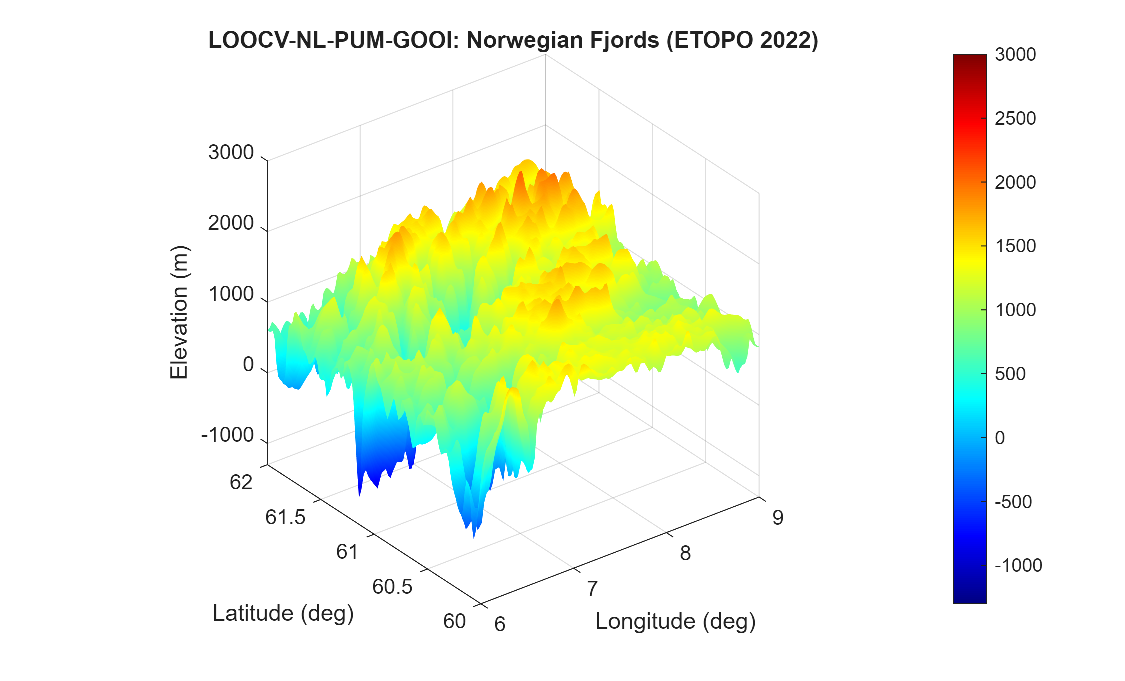}
\end{minipage}
\caption{Top-down view (left) and 3D view (right) of the Norwegian Fjords
elevation comparative reconstructions (ETOPO~2022, Sognefjord region).
Rows from top: raw data, Linear PUM, NL-PUM, LOOCV-NL-PUM-GOOI.}
\label{fig:etopo}
\end{figure}

\Cref{fig:etopo} (left) displays top-down views, while \cref{fig:etopo}
(right) shows the corresponding 3D surfaces. It can be observed that the
Linear PUM fails to reconstruct the surface accurately and overshoots above
$2500$~m and creates phantom valleys near $-2000$~m, which confirms the
occurrence of the Gibbs phenomenon across the sharp land-sea elevation
changes. However, although NL-PUM eliminates the oscillations, it
over-smooths the fjord channels, reducing their depth and blurring sharp
transitions. LOOCV-NL-PUM-GOOI yields the most accurate reconstruction:
the fjord channels near $61.2^{\circ}$N are sharply resolved, mountain
peak heights are well preserved, and no more Gibbs phenomenon is visible.
This visual improvement is consistent with the quantitative outperformance
observed in the test functions above: the adaptive radius contracts near
sharp elevation gradients, preventing the local RBF interpolant from mixing
land and sea data, while remaining at its initialized value in smooth
highland or open-sea regions.

\clearpage
\newpage

\section{Conclusion}
\label{sec:conclusion}

In this work, we have introduced LOOCV-NL-PUM-GOOI, an adaptive extension
of the NL-PUM method \cite{ramon25} that incorporates two novelties: a
local smoothness indicator $\sigma_j$ based on the normalized mean absolute
residual of a degree-$2$ polynomial fit, and an exponential radius
shrinkage mechanism that automatically contracts patches near
discontinuities without requiring any prior knowledge of the interface
geometry. In smooth regions, the method reduces to the LOOCV-GOOI parameter
selection strategy of \cite{cav25b} applied locally within the PU
framework, with no additional modifications from the discontinuity-handling
components.

The numerical results lead to two clear conclusions. First, adapting the
shape parameter alone without modifying the patch radius is insufficient
near discontinuities, as confirmed by the LOOCV-NL-PUM results. Second,
the proposed radius shrinkage mechanism consistently reduces the root mean
square error across all discontinuous test cases while remaining
computationally comparable to NL-PUM. The application to Norwegian Fjords
elevation data confirms that these gains extend to real problems with
irregular point distributions and multiple sharp transitions. The method
does present limitations. Strong radius shrinkage near compact interior
interfaces can elevate the maximum pointwise error at isolated locations, a
known trade-off of radius-based strategies. A natural direction for future
work is the design of interpolation schemes that simultaneously handle jump
discontinuities and spatially varying point density, addressing the two
principal sources of difficulty in real-world scattered data problems and
extending the framework developed here in a unified direction.

\section*{Reproducibility}

The MATLAB codes necessary to reproduce all numerical results, figures, and 
tables presented in this paper are available in the repository 
\url{https://github.com/ahaider97/NL-PUM-LOOCV-GOOI}. 
The repository contains the implementations of all four methods compared 
in this work, along with the scripts used to generate each figure 
and table. The real data application uses 
the ETOPO 2022 global relief model \cite{etopo2022}, which must be 
downloaded separately from \url{https://doi.org/10.25921/fd45-gt74}; instructions are provided in the repository README.
\section*{Acknowledgments}
\begin{sloppypar}
The work of A.H.\ and R.C.\ has been supported by the INdAM Research group
GNCS, the GFI 2025 Project, and the 2024 Project ``Numerical Analysis and
Modelling'' funded by the Department of Mathematics ``Giuseppe Peano'' of
the University of Torino. This research has been accomplished within the
RITA ``Research ITalian network on Approximation,'' the UMI Group TAA
``Approximation Theory and Applications,'' and the SIMAI Activity Group
ANA\&A ``Numerical and Analytical Approximation of Data and Functions with
Applications.'' The work of A.H.\ and D.F.Y.\ has been supported by grant
PID2023-146836NB-I00, funded by MICIU/AEI/10.13039/501100011033 and by
ERDF/EU. The work of D.F.Y.\ and J.R.-A.\ has been supported by GVA project
CIAICO/2024/089.
\end{sloppypar}
\bibliographystyle{siamplain}
\bibliography{references}

\end{document}